\documentclass[10pt]{amsart}
\pdfoutput=1
\usepackage{amssymb}
\usepackage{booktabs}
\usepackage{mathtools}
\usepackage{anyfontsize}
\usepackage{multicol}
\usepackage{tikz-cd}
\usepackage{adjustbox}
\usepackage{array}
\usepackage{rotating}
\usepackage{booktabs}
\newenvironment{amssidewaystable}
  {\begin{sidewaystable}\vspace*{.5\textwidth}\begin{minipage}{\textheight}\centering}
  {\end{minipage}\end{sidewaystable}}

\usepackage[
    backend=biber
    ,style=alphabetic
    ,sorting=nyt
]{biblatex}
\usepackage[colorlinks,
    citecolor=green,
    filecolor=black,
    linkcolor=blue,
    urlcolor=black,
]{hyperref}
\hypersetup{linktocpage}
\usepackage[
    numbered,
    open=true,
]{bookmark}

\newcolumntype{R}[2]{%
    >{\adjustbox{angle=#1,lap=\width-(#2)}\bgroup}%
    l%
    <{\egroup}%
}

\theoremstyle{plain}
\newtheorem{theorem}{Theorem}
\newtheorem{corollary}[theorem]{Corollary}
\newtheorem{proposition}[theorem]{Proposition}
\newtheorem{lemma}[theorem]{Lemma}

\theoremstyle{definition}
\newtheorem{definition}[theorem]{Definition}

\newtheorem{remark}[theorem]{Remark}
\newtheorem{example}[theorem]{Example}

\numberwithin{theorem}{section}

\usepackage{tikz}
\usetikzlibrary{positioning}
\usetikzlibrary{decorations.markings}
\usetikzlibrary{decorations.pathreplacing}
\usetikzlibrary{backgrounds}
\usetikzlibrary{calc}

\usepackage{orcidlink}

\newcommand{\Sgrp}{\mathsf{Sgrp}}
\newcommand{\Grp}{\mathsf{Grp}}

\newcommand{\N}{\mathbb{N}}
\newcommand{\Z}{\mathbb{Z}}

\renewcommand{\geq}{\geqslant}
\renewcommand{\leq}{\leqslant}
\renewcommand{\l}{\ell}
\renewcommand{\phi}{\varphi}
\renewcommand{\epsilon}{\varepsilon}
\newcommand{\id}{\operatorname{id}}
\newcommand{\tensor}{\mathrel{\otimes}}

\newcommand{\op}{\operatorname{op}}
\newcommand{\lcm}{\operatorname{lcm}}
\newcommand{\rank}{\operatorname{rank}}

\newcommand{\Tor}{\operatorname{Tor}}

\newcommand{\bbone}{\text{\usefont{U}{bbold}{m}{n}1}}
\MakeRobust{\bbone}
\newcommand{\one}{^\bbone}

\newcommand{\im}{\operatorname{im}}

\begin{document}

\title[Homotopy Types of Small Semigroups]{Homotopy Types of Small Semigroups}
\author{Dennis Sweeney}

\begin{abstract}
    We present a software package
    for quickly calculating the integral
    homology of finite semigroups and monoids by a novel approach
    of exploiting structure in projective resolutions.
    We describe the usage of this package
    on semigroups and monoids of small orders.
    From the results of these calculations,
    we present a wealth of counterexamples in finite semigroup theory:
    we give a finite aperiodic semigroup with
    large torsion in its homology,
    two finite semigroups with certain Moore spaces
    as classifying spaces,
    and a finite semigroup and with
    nontrivial rational homology in infinitely many dimensions,
    refuting conjectures of William Nico.
    We further show that the set of homotopy types of
    classifying spaces of finite semigroups is closed under suspension.
\end{abstract}

\maketitle
\tableofcontents

\subsection*{Acknowledgements}

The author would like to thank Sanjeevi Krishnan for many helpful conversations
and suggestions.

\section{Introduction}\label{sec:intro}
The classifying space $BG$ of a (discrete) group $G$ is
an Eilenberg-MacLane space of type $K(G,1)$: it satisfies
$\pi_1(BG) \cong G$, and $\pi_k(BG)=0$ for $k \neq 1$.
However, \cite{McDuff79} showed that
the classifying space $BM$ of a monoid $M$
has no such homotopical restriction: $BM$
can attain the homotopy type of any connected CW complex.

Finite semigroups and monoids, and especially their actions
on a finite set of states, have long been
studied in automaton theory \cite[]{KrohnRhodes,ArbibKrohnRhodes},
and much is known about the properties and structure of finite semigroups
\cite[]{GreensRelations,CliffordPrestonVol1,qTheory}.
Fiedorowicz conjectured \cite[p. 323]{Zig84}
that any finite simply connected
CW complex is homotopy equivalent to the
classifying space of some finite monoid,
but to date, the classifying spaces of the finite monoids considered
have generated relatively few unique homotopy types
(though this is growing \cite[]{realizingMooore}).
The classifying spaces of finite groups are well-studied \cite[]{brown_coho},
and \cite{Zig84,Zig02}
notes that a $2$-by-$2$ rectangular band semigroup
(with an adjoined $1$)
has a classifying space homotopy equivalent
to the $2$-sphere $\mathbb{S}^2$,
while \cite{steinberg2024} computes more generally the classifying spaces of completely simple semigroups with techniques from \cite{GraySteinberg2022}.
Although \cite{adjoint_derived} and \cite{nico69} consider certain reductions
to smaller monoids, and \cite{nicoHomologicalDimension,nicoImproved}
gives bounds on cohomology of finite \textit{regular} semigroups
in high dimensions,
the exact homotopy types of classifying spaces of individual
general finite semigroups nonetheless remain elusive.

To begin to rectify this, we uniformly investigate
semigroups of small orders.
We use the \texttt{minion} solver \cite[]{minion_paper, minion} to generate
finite semigroups and monoids,
and then we use a new
\texttt{fast\_semigroup\_homology}
Python package developed by the author \cite[]{fast_semigroup_homology}
to compute homology using projective resolutions.
Our novel methods of constructing
small projective resolutions,
splitting resolutions into multiple branches,
and caching pieces of resolutions
offered substantial performance improvements compared
to computations involving the bar resolution
(see Table~\ref{tab:benchmarks}),
which allowed many more homology groups to be computed.

Using the \texttt{fast\_semigroup\_homology} package, we compute homology groups
$H_0(S),\dots, H_8(S)$ for all semigroups of order at most 8,
and then we reduce to a restricted
class of monoids (those with no one-sided ideals that are monoids)
and compute several homology groups of such monoids of order at most 11.
We also search through a further restricted class of monoids
(those with no one-sided ideals that are monoids and
 at most 4 non-unit elements outside the minimal ideal)
and continue computing homology for such monoids of order at most 16.

For finite nontrivial groups $G$, the reduced (co)homology of $BG$
is annihilated by multiplication by $|G|$
\cite[Corollary 10.2]{brown_coho}, and $BG$
has nontrivial (co)homology groups in infinitely many dimensions
\cite[]{swan59}. Refuting semigroup generalizations of these facts,
our computations reveal a wealth of counterexamples
in finite semigroup theory, summarized in Table~\ref{tab:false_conjectures}.

\begin{table}[ht]
    \caption{Claims that hold for finite groups
            but do not hold for finite semigroups or monoids,
            with references to the specific counterexamples
            presented in this paper.}
    \label{tab:false_conjectures}
\begin{tabular}{ll}
    \toprule
    \textbf{False Conjecture on Finite $S$}
    & \textbf{Counterexample}
    \\ \midrule
    \begin{tabular}{@{}l@{}}
        If $BS$ is contractible,\\
        then $S$ has a left or right zero.
    \end{tabular}
    &
    \begin{tabular}{@{}l@{}}
        \cite{adjoint_derived} has Ex. 3\\
        with $\geq 9$ elements. \\
        Our \ref{ex:no_zero_contractible} has $|S|=5$. \\
    \end{tabular}
    \\ \midrule
    \begin{tabular}{@{}l@{}}
        Subgroups of $S$ determine $H^i(S)$ for large $i$; \\
        $S$ aperiodic $\Rightarrow |H^i(S)| < \infty$ for large $i$. \\
        Conjectured on p.598 of \cite{nico69}.
    \end{tabular}
    & \ref{ex:0ZZZ} has $H_i(S)=\Z$ for $i\geq 2$.
    \\ \midrule
    \begin{tabular}{@{}l@{}}
        The minimum $\#$ of generators for $H_i(S)$
        \\
        $\leq \text{polynomial}(i)$. \\
        Holds for groups by \cite{evens_fingen_coho}.
    \end{tabular}
    & \ref{ex:2pow} has $H_i(S)=\Z^{2^{i-2}}$ for $i\geq 2$.
    \\ \midrule
    \begin{tabular}{@{}l@{}}
        If $H_i(S)$ has $p$-torsion for some $i$,
        then \\
        $H_i(S)$ has $p$-torsion for infinitely many $i$.\\
        Holds for groups by \cite[Cor. 1]{swan59}.
    \end{tabular}
    &
    \begin{tabular}{@{}l@{}}
        \ref{ex:Moore_C2_3} has $BS \simeq \mathsf{M}(C_2, 3)$. \\
        \ref{ex:Moore_C2_2} has $BS \simeq \mathsf{M}(C_2, 2)$.
    \end{tabular}
    \\ \midrule
    \begin{tabular}{@{}l@{}}
        If $H_i(S)$ has $p$-torsion, \\
        then $p\leq |S|$, and $S$ has a $p$-subgroup. \\
        Holds for groups by \cite[Cor. 10.2]{brown_coho}.
    \end{tabular}
    &
    \begin{tabular}{@{}l@{}}
        \ref{ex:Moore_C2_3} has $BS \simeq \mathsf{M}(C_2, 3)$. \\
        \ref{ex:1494640} has $H_6(S) = \Z^9 \times C_{1494640}$. \\
    \end{tabular}
    \\ \midrule
    \begin{tabular}{@{}l@{}}
        Nontrivial cohomology $H^*(S)$ has\\
        some nontrivial cup product. \\
        Holds for groups by \cite{evens_fingen_coho}.
    \end{tabular}
    & \ref{ex:sus_C2} has $BS \simeq \Sigma(BC_2)$.
    \\ \midrule
    \begin{tabular}{@{}l@{}}
        The ring $H^*(S)$ is finitely generated. \\
        Holds for groups by \cite{evens_fingen_coho}.
    \end{tabular}
    &
   \begin{tabular}{@{}l@{}}
        \ref{ex:2pow} has $H_i(S) = \Z^{2^{i-2}}$ for $i \geq 2$. \\
        \ref{ex:sus_C2} has $BS \simeq \Sigma(BC_2)$. \\
    \end{tabular}
    \\ \bottomrule
\end{tabular}
\end{table}

The homology of the 6-element monoid of Example~\ref{ex:sus_C2}
was computed to be that of $\Sigma(BC_2)$,
the suspension of the classifying space of the cyclic group of order $2$.
In fact, the classifying space of Example~\ref{ex:sus_C2} is
homotopy equivalent $\Sigma(BC_2)$, and we generalize to prove the following:
\begin{theorem}\label{thm:join_semigroup}
    For any monoid $M$ and any set $Y$, there is a monoid $J^Y(M)$
    of order $(1+|Y|)|M|$ with $B\!\left(J^Y(M)\right)\simeq BM * Y^{\mathrm{discrete}}$, the topological join of $BM$ with the discrete set $Y$.
\end{theorem}
This has the following corollaries:
\begin{corollary}\label{cor:wedge_of_suspensions}
    For any nonempty set $Y$ and $y_0 \in Y$,
    and any monoid $M$,
    the monoid $J^Y(M)$ has classifying space
    \[
        B\!\left(J^Y(M)\right) \simeq \bigvee_{y \in Y \setminus \{y_0\}} \Sigma(BM),
    \]
    where $\Sigma$ denotes suspension.
\end{corollary}
\begin{corollary}\label{cor:suspension}
    For any monoid $M$, there is a monoid
    $M' \coloneq J^{\{1,2\}}(M)$ with $|M'|=3 \cdot |M|$
    and $BM' \simeq \Sigma(BM)$,
    so the set of homotopy types of classifying spaces
    of finite semigroups is closed under suspension.
\end{corollary}
\begin{corollary}\label{cor:wedge_of_spheres}
    For any $n \geq 2$ and $r \geq 0$, there is a finite
    monoid $M$ with classifying space
    $BM \simeq \bigvee_{k=1}^{r} \mathbb{S}^n$,
    a wedge of $r$ copies of the $n$-sphere.
\end{corollary}
Section~\ref{sec:background} below recalls facts about semigroups, monoids,
simplicial sets, and homological algebra, as applicable to our computations.
Section~\ref{sec:fast} discusses the design
of the \texttt{fast\_semigroup\_homology} Python package \cite[]{fast_semigroup_homology},
and then Section~\ref{sec:finding} discusses the design of the \texttt{semisearch}
program \cite[]{semisearch} for generating semigroups and monoids.
Section~\ref{sec:results} presents the results of the computations,
and Section~\ref{sec:suspensions} gives proofs of
Theorem~\ref{thm:join_semigroup} and its corollaries.

\section{Theoretical Background}\label{sec:background}
\subsection{Semigroups and Monoids}\label{sub:background_semigroups}

For a general background in semigroup theory, see
\cite{HowieSemigroups}, \cite{CliffordPrestonVol1},
Section 1.1 of \cite{HindmanStrauss},
or Section 1 of \cite{repFiniteMon}.

A semigroup $S$ is a set equipped with an associative binary operation.
A monoid $M$ is a semigroup that contains an identity element, denoted $1$.
Given a semigroup $S$, we will write $S\one$ for
the monoid formed by adjoining an identity element to $S$,
even if $S$ already has an identity, so $S\one = S \sqcup \{1\}$
with $1x=x1=1$ for all $x\in S\one$. A \textit{subgroup} of $S$
is a subsemigroup isomorphic to a group. An \textit{aperiodic} semigroup $S$
is a semigroup such that each $x\in S$ satisfies $x^{n+1}=x^n$ for some $n\in\N$.
A finite semigroup is aperiodic if and only if it has only trivial subgroups.
For a semigroup
$S$, we will write $S^{\op}$ for the semigroup with the same underlying set as $S$
but with the dual operation $x \cdot y \coloneq yx$. An element $z \in S$
is called a \textit{left zero} (resp. \textit{right zero}) if $zx=z$ (resp. $xz=z$)
for all $x \in S$.

A \textit{left ideal} (resp. \textit{right ideal}) of a semigroup $S$ is
a nonempty subset $A \subseteq S$ satisfying $SA \subseteq A$ (resp. $AS \subseteq A$).
An \textit{ideal} of $S$ is a subset that is both a left ideal
and a right ideal.
Finite intersections
of ideals are ideals, so if $S$ is finite and nonempty, then the intersection
of all ideals of $S$ is an ideal. This \textit{minimal ideal} (or ``kernel'' on p.67 of \cite{CliffordPrestonVol1})
of a finite nonempty semigroup $S$ is denoted $K(S)$,
and $K(S)$ is a \textit{simple} semigroup in the sense that it has no proper ideals.

Per \cite{sus28} and \cite{rees1940} (see also Section 3.2 of \cite{HowieSemigroups} or Theorem~1.64 of \cite{HindmanStrauss}),
the classification of finite simple (or more generally, \textit{completely simple})
semigroups is well-understood, up to the classification of groups:
each finite simple semigroup
is isomorphic to a \textit{Rees matrix semigroup}, denoted
$
    \mathcal{M}(H;I,J;C_{(-,-)}),
$
where $H$ is some finite group, $I$ and $J$ are finite sets,
and $C_{(-,-)}\colon J \times I \to H$ is any function,
called the \textit{sandwich matrix}. This semigroup
has the underlying set $I\times H \times J$ with
operation
$
    (i, h, j)(i', h', j') \coloneq (i, h C_{(j,i')}h', j').
$
Each $\{i\} \times H \times \{j\}$ is then a subgroup isomorphic to $H$,
with identity $(i, C_{(j,i)}^{-1}, j)$.
Furthermore, every finite simple semigroup
is isomorphic to a Rees matrix semigroup with a \textit{normalized}
sandwich matrix, meaning that there is some $i_0 \in I$
and $j_0 \in J$ such that for all $i\in I$ and $j\in J$ we have
$C_{(j_0, i)} = C_{(j, i_0)} = 1_H$.
If $H$ is trivial then $\mathcal{M}(H;I,J;C_{(-,-)})$
is called an $|I|$-by-$|J|$ \textit{rectangular band}, while if $C_{(-,-)}$
is a constant function then $\mathcal{M}(H;I,J;C_{(-,-)})$
is a direct product of a rectangular band and a group.
More generally, a Rees matrix semigroup can be understood
as a ``twisted product'' of a rectangular band and a group.

For any finite nonempty semigroup $S$, if we write
$K(S) = \mathcal{M}(H;I,J;C_{(-,-)})$,
then each set $I \times H \times \{j\}$
is a minimal left ideal of $S$,
and each set $\{i\} \times H \times J$ is a minimal right ideal;
see the statement of Theorem~1.64 in \cite{HindmanStrauss}.
Thus, for each $(i, h, j) \in K(S)$ and $x \in S$,
we have $(i, h, j)x = (i, h', j')$ for some $h' \in H$ and $j' \in J$,
and likewise $x(i, h, j) = (i'', h'', j)$ for some $i'' \in I$ and $h'' \in H$.

The \textit{group completion} of a semigroup $S$
is a group $GS$ with a semigroup homomorphism $\gamma\colon S \to GS$
such that any semigroup homomorphism $S\to \Gamma$ to a group
$\Gamma$ factors uniquely as $S\xrightarrow{\gamma} GS \to \Gamma$.
The group completion
is unique up to isomorphism, and the functor $G\colon \Sgrp \to \Grp$
from the category of semigroups to the category of groups
is the left adjoint of the forgetful functor $\Grp \to \Sgrp$
from groups to semigroups. The group $GS$ can also
be called the \textit{maximal homomorphic group image} of $S$
\cite[]{HM_group_image}, the \textit{enveloping group} of $S$
\cite[]{mcduff_segal_group_completion},
or when $S$ is a commutative monoid, the \textit{Grothendieck group} of $S$.
The group completion of a finite semigroup is readily computed:
\begin{lemma}[\cite{larisse_group_image}, Lemme II.3; \cite{HM_group_image}, p. 848]
    \label{lem:group_completion}
    For any finite nonempty semigroup $S$, identifying the minimal
    ideal $K(S)$ with some Rees matrix semigroup $\mathcal{M}(H;I,J;C_{(-,-)})$
    with normalized sandwich matrix as described above, we have
    \begin{align*}
        GS = G(K(S))
        &= H/\langle\langle C_{(j,i)} \mid j\in J, i \in I\rangle\rangle,
        \\&= H/\langle\langle C_{(j,i)} \mid j\in J\setminus \{j_0\}, i \in I\setminus \{i_0\}\rangle\rangle
    \end{align*}
    where $\langle\langle {\cdots} \rangle\rangle$ denotes the normal closure
    $N\trianglelefteq H$
    of a subset of $H$, and where the group completion map $S \to GS$
    sends each $x \in S$ to $[h_x] \in H/N$
    such that $(i_0, 1_H, j_0)x(i_0, 1_H, j_0) = (i_0, h_x, j_0)$.
\end{lemma}
Note in particular that the group completion of a finite semigroup is finite,
and $G(S\one) = GS$. The group completion of a monoid is the same as its group completion
as a semigroup.

For a monoid $M$, the \textit{group of units}
is the set of $x \in M$ such that some $y,z \in S$ satisfy
$yx=xz=1$. The group of units is the unique maximal subgroup of $M$ containing $1$
\cite[Theorem 1.10]{CliffordPrestonVol1}. If $M$ is a finite monoid,
the set of non-units is an ideal of $M$ \cite[Exercise A.2.2]{qTheory}.

\subsection{Classifying Spaces}\label{sub:classifying}

The classifying space of a monoid is the realization of a simplicial set,
and the classifying space of a general semigroup is more naturally described as a Delta set.
For an introduction to Delta sets and simplicial sets,
see \cite{ssets}, or see
Section 2.1 and the \textit{Simplicial CW Structures} appendix of \cite{HatcherAT}.
Simplicial sets are discussed more thoroughly in
Section 8.1 of \cite{weibelHA} and in \cite{maySimplicial}.
For an overview of classifying spaces of categories,
refer to \cite{segal68} or Section IV.3 of \cite{kbook}.

Recall that a \textit{Delta set} $X$
is a sequence of sets
$X^0, X^1, X^2, \dots$
equipped with face maps $d^i_j \colon X^i \to X^{i-1}$ satisfying certain identities,
and a \textit{simplicial set} is a Delta set additionally
equipped with degeneracy maps $s^i_j\colon X^i \to X^{i+1}$
satisfying additional identities.

The \textit{nerve} $NM$ of a monoid $M$
is the same simplicial set as its nerve when viewed
as a category with one object, and thus has an $n$-simplex for each
list of $n$ elements of $M$. We will write such an $n$-simplex as $[m_1|\cdots|m_n]$.
The face maps of $NM$
remove the vertical bars in the notation to replace $m_i |m_{i+1}$
with the product $m_im_{i+1}$ (or delete $m_1$ or $m_n$),
and the degeneracy maps insert the identity $1 \in M$,
so any $n$-simplex that includes $1$ is degenerate.
We then define the classifying space $BM := |NM|$ to be the
geometric realization of the nerve,
working in the category of compactly generated weak Hausdorff spaces.
Recall that natural transformations between functors
descend to homotopies between classifying spaces \cite[Prop. 2.1]{segal68}.

If $S$ is a semigroup instead of a monoid, we could similarly define its
classifying space $B_\Sgrp S$
as the realization of the Delta set with similar $n$-simplices $[x_1|\cdots|x_n]$
and face maps, but without degeneracy maps, but this is identical to $B(S\one)$.
Furthermore, if $S$ is a monoid, then there is a homotopy equivalence
$BS \simeq B_\Sgrp S$, since
these are the thick and standard realizations of the same simplicial set.
Hence, there is no essential difference
between results here on semigroups and on monoids,
and we will simply write $BS \coloneq B_\Sgrp(S) \cong B(S\one)$.

\subsection{Homology as Tor}
Since \cite{McDuff79} showed that connected CW homotopy types have discrete
algebraic representations as monoids, it is desirable
to find correspondences between algebraic properties
of monoids and topological properties of their classifying spaces.
The first well-known result in this direction is that
the fundamental group of $BM$ is $GM$ \cite[Lemma 1]{McDuff79}.
Another class of results, described for example in \cite{Nunes1995},
is that some standard algebraically-defined homology
theories for $M$ agree with the homology of $BM$.

For a monoid $M$, we define the \textit{monoid ring} $\Z M$ to be the (unital) ring
of finite formal $\Z$-linear combinations of $M$ elements.
Given a left $\Z M$-module $D$, we define the
homology of $M$ with coefficients in $D$:
\[
    H_i(M; D) \coloneq \Tor^{\Z M}_i(\Z, D).
\]
Above, $\Z$ is a right $\Z M$-module with trivial $\Z M$ action.
This is identical to the definition given in Section X.5 of \cite{homology},
with the exception that we use a left module of coefficients,
as is done for groups in Chapter 6 of \cite{weibelHA}
and Section 9.4 of \cite{rotmanHA}.
We will most often consider the case $D=\Z$ and write
\[
    H_i(M) \coloneq H_i(M; \Z) = \Tor^{\Z M}_i(\Z, \Z).
\]
We now have the following:
\begin{proposition}[From Proposition 0.4.2 in \cite{Nunes1995}]
    \label{localCoefficients}
    For a monoid $M$, and any $GM$-module $D$, we have
    \[
        H_i(BM; D) \cong H_i(M; D),
    \]
    where the left side is interpreted as homology
    with local coefficients,
    and $D$ is interpreted as a $\Z M$-module via the group completion
    map $\gamma \colon M \to GM$, so $m\cdot v \coloneq \gamma(m)v$
    for $m \in M$ and $v \in D$.
\end{proposition}
In particular, we have the following, which is well-known (e.g. by \cite{Zig02}):
\begin{corollary}\label{cor:opposite_equivalent}
    For any monoid $M$, we have isomorphic integral homology
    \[
        H_i(M) \cong H_i(BM) \cong H_i(BM^{\op}) \cong H_i(M^{\op}).
    \]
\end{corollary}
\begin{proof}
    For the first and third isomorphisms, take $D=\Z$ in the previous proposition.
    For the middle, there is an obvious homeomorphism $BM \cong BM^{\op}$.
\end{proof}
\begin{corollary}
    For any semigroup $S$, the integral semigroup homology,
    defined as $H_i(S) \coloneq H_i(S\one)$,
    is isomorphic to $H_i(B_{\Sgrp} S)$.
\end{corollary}
\begin{proof}
    Since $B(S\one) \cong B_{\Sgrp}S$,
    we have $H_i(S\one) \cong H_i(BS\one) \cong H_i(B_{\Sgrp} S)$.
\end{proof}
\begin{remark}
    As pointed out by the proof of Proposition 4.4 of \cite{Zig84},
    the universal cover $\tilde{BM}$ of the classifying space $BM$ of a monoid
    has homology $H_i(\tilde{BM}; \Z) \cong \Tor^{\Z M}_i(\Z, \Z GM)$.
    By the Hurewicz theorem, we have
    \[
        \pi_2(BM) \cong \pi_2(\tilde{BM}) \cong H_2(\tilde{BM}) \cong \Tor^{\Z M}_i(\Z, \Z GM),
    \]
    so the computations described in Section~\ref{sec:fast}
    can be easily adapted to compute the second homotopy group
    of the classifying space of a finite semigroup or monoid.
\end{remark}

\subsection{Reduction to Smaller Semigroups}
\cite{adjoint_derived}
give the following way to reduce homology computations
of large semigroups to homology computations of certain subsemigroups
(see also Lemma 2.1 of \cite{nico69}):
\begin{lemma}[Theorem 2 of \cite{adjoint_derived}]
    \label{lem:eSeHomology}
    Let $S$ be a semigroup and let $M$
    be a left (right) ideal of $S$ that is a monoid.
    Let $D$ be any left (right) $S$-module (defined as a $\Z S\one$-module).
    Let $e$ be the identity element of $M$. Then for
    all $n \geq 0$, we have
    \[
        H^n(S, D) \cong H^n(M, D) = H^n(S, eD) = H^n(M, eD).
    \]
\end{lemma}
Based on this, with the note that the first isomorphism above
is induced by the map $S\to M$ given by $x\mapsto exe$,
we see that this semigroup homomorphism induces
a function $BS \to BM$ that induces isomorphisms
on cohomology with all local coefficient systems
(i.e. $\Z G(S)$-modules) $D$, so this function
is a homotopy equivalence $BS \simeq BM$.
The following is an original direct proof of the same fact:
\begin{theorem}
    \label{thm:eSeEquivalent}
    If $S$ is a semigroup with some idempotent $e$ such that $eSe=Se$
    or $eSe=eS$, then $BS \simeq B(eSe)$.
\end{theorem}
\begin{proof}
    By replacing $S$ with $S^{\op}$ as necessary,
    assume without loss of generality that $eSe=eS$.
    Write $M\coloneq eS$,
    and note that $e$ is an identity for $M$,
    so $exe=ex$ for all $x\in S$. Define $r\colon S \to M$
    by $r(x)=ex$. The map $r$ is a semigroup
    homomorphism because $exey=exy$. Write $\iota\colon M\to S$
    for the inclusion, so then $r\circ \iota = \id_M$.
    We will show that $\iota$ and $r$ define a homotopy equivalence.
    Since $B_{\Sgrp} S \cong B S\one$ is a natural homeomorphism,
    it suffices to show that $r\one \colon S\one \to M\one$
    and $\iota\one \colon M\one \to S\one$
    define a homotopy equivalence $B r\one \colon BS\one \to B M\one$
    and $B \iota\one \colon BM\one \to BS\one$.

    Interpreting $M\one$ and $S\one$ as categories
    each with a single object, we provide natural transformations as homotopies.
    The first composition $r\one \circ \iota\one = (r\circ \iota)\one = \id_M\one = \id_{M\one}$ is the identity functor already, so no homotopy is needed.
    For the second, define a natural transformation
    $\eta$ from the identity functor on $S\one$
    to $\iota\one \circ r\one$ by taking the unique component
    $\eta_* = e$. The naturality square for each $x \in S\one$ is:
    \[
    \begin{tikzcd}
        * \ar[r, "\eta_*=e"]
        \ar[d, "\id_{S\one}x=x" left]
        &*\ar[d, "\iota\one r\one x"]
        \\*\ar[r, "\eta_*=e" below]
        &*
    \end{tikzcd}
    \]
    This commutes because $e1=1e$ for $x=1$
    and $ex=exe$ for $x \in S$.
\end{proof}

The hypothesis above holds exactly when $S$ has a one-sided ideal that is a monoid
with identity $e$, or equivalently when there is
an idempotent $e$ such that $exe=ex$ (or $exe=xe$) for all $x \in S$.
We have the following major consequence, the cohomology version
of which was also pointed out as Corollary 4 of \cite{adjoint_derived}:
\begin{corollary}
    \label{cor:kThin}
    If $S$ is a semigroup with minimal ideal $K(S)\cong \mathcal{M}(H; I, J; C_{(-,-)})$
    in which $|I|=1$ or $|J|=1$,
    then $BS \simeq BH$.
\end{corollary}
\begin{proof}
    Without loss of generality write $K(S) = \mathcal{M}(H; I, J, C_{(-, -)})$
    with $J = \{j_0\}$.
    Then $S$ has a minimal right ideal
    \[
        (i_0, 1_H, j_0)S = \{i_0\} \times H \times J = \{i_0\} \times H \times \{j_0\},
    \]
    so $eSe=eS$, where $e=(i_0, 1_H, j_0)$.
    But $eSe\cong H$, so $BS \simeq B(eSe) \cong BH$.
\end{proof}
This applies for large classes of finite semigroups, as is well-known:
\begin{corollary}\label{cor:zero_contract}
    If $S$ is a semigroup with a left or right zero element,
    then $BS$ is contractible.
\end{corollary}
\begin{proof}
    If $z$ is a left zero, then $zS=\{z\}=zSz$,
    so $BS \simeq B(\{z\})$. The classifying space
    of the trivial monoid is contractible because it has
    only one nondegenerate simplex, the $0$-simplex ``[]''.
    The argument for right zeros is identical.
\end{proof}
\begin{corollary}
    If $S$ is a finite semigroup,
    and $S$ satisfies some equation
    \[
        x_{\lambda_1}\cdots x_{\lambda_\l} = x_{\lambda'_1}\cdots x_{\lambda'_{\l'}}
    \]
    for all substitutions of semigroup elements
    into the variables $x_1, x_2, x_3, \dots$,
    and either $\lambda_1\neq \lambda'_1$ or $\lambda_\l \neq \lambda'_{\l'}$,
    then $BS \simeq BH$, where $H$ is the finite group defining $K(S)$
    as above.
\end{corollary}
\begin{proof}
    Replacing $S$ with $S^{\op}$ as necessary,
    assume $\lambda_1 \neq \lambda_1'$.
    Write $K(S)$ as $\mathcal{M}(H; I,J; C_{(-,-)})$.
    For any $(i, h, j), (i',h',j') \in K(S)$,
    substitute $x_{\lambda_1} = (i,h,j)$
    and $x_{\lambda'_1}=(i',h',j')$.
    Then $x_{\lambda_1}\cdots x_{\lambda_\l}$
    evaluates to some element of $\{i\} \times H \times J$,
    while $x_{\lambda'_1} \cdots x_{\lambda'_{\l'}}$
    evaluates to some element of $\{i'\} \times H \times J$,
    but these must be equal, so $i=i'$. This shows that $|I|=1$.
\end{proof}
For example, \textit{right-commutative} finite semigroups (which satisfy the equation $x_1x_2x_3=x_1x_3x_2$)
have classifying space homotopy equivalent to that of some group.

Note also that finite \textit{inverse} semigroups have commuting idempotents
\cite[Theorem 1.17]{CliffordPrestonVol1},
so the minimal ideal of a finite inverse semigroup must have $|I|=|J|=1$,
and so Theorem~\ref{thm:eSeEquivalent} also applies to finite inverse semigroups.

\section{The \texttt{fast\_semigroup\_homology} Package}\label{sec:fast}
The \texttt{fast\_semigroup\_homology} package used to compute the homology
of finite semigroups is available at \cite{fast_semigroup_homology}.
This section describes its design.

Given a finite semigroup, we can apply Theorem~\ref{thm:eSeEquivalent}
to reduce to a smaller monoid or else adjoin an identity
to get a monoid with the same classifying space,
so for the rest of this section
we will consider a finite monoid $M$ defined by a multiplication table.

To compute the integral homology $H_i(M)$,
we build a projective resolution $\Z \leftarrow P_0 \leftarrow P_1 \leftarrow \cdots$
of the trivial left $\Z M$-module $\Z$, then
delete the augmentation ($\Z \leftarrow$),
tensor the complex with the trivial right $\Z M$-module $\Z$,
and take the homology of the resulting chain complex:
\[
    H_i(BM) = H_i(M) = \Tor_i^{\Z M}(\Z, \Z)
    = H_i\!\left(\Z \tensor_{\Z M} P_0 \leftarrow \Z \tensor_{\Z M} P_1 \leftarrow \cdots \right).
\]

The major challenge of the computation
is to find sufficiently small projective resolutions.
To build these resolutions, we will use projective modules
of the form $P_n = \Z M e_1 \oplus \cdots \oplus \Z M e_{\l_n}$,
where $e_i$ are idempotents of $S$. Each
$\Z M e_i$ is projective because the monoid ring $\Z M$
decomposes as a direct sum $\Z M = \Z M e_i \oplus \Z M (1 - e_i)$
of left $\Z M$-modules, and direct sums of projective modules are projective.

The maps $P_{n} \to P_{n-1}$ of our projective resolutions
will be maps
\[
    d \colon \bigoplus_{j=1}^{m} \Z M e_j \to \bigoplus_{i=1}^{m'} \Z M f_i
    \quad
    \text{defined by}
    \quad
    p_i(d(\iota_j(x))) = x a_j^i,
\]
where $e_1, \dots e_m, f_1, \dots, f_{m'}$ are idempotents of $M$,
$p_i$ is projection onto summand $i$,  $\iota_j$ is inclusion
of summand $j$, and $(a_j^i)$ is a matrix of elements
of $\Z M$. For $d$ to have the listed codomain,
the elements $a_j^i$ must be chosen to satisfy
$e_j a_j^i \in \Z M f_i$,
so $xe_j a_j^i \in x (\Z M f_j) \subseteq \Z M f_j$
for any $x e_j \in \Z M e_j$. The left $\Z M$ linearity
of the map $d$ follows immediately from the associativity
of the monoid operation.

\begin{example}[p.34 of \cite{Zig02}]
    Taking the $2$-by-$2$ rectangular band
    semigroup $S=\{x_{00}, x_{01}, x_{10}, x_{11}\}$
    with operation $x_{ij}x_{i'j'} = x_{ij'}$,
    we find the following projective resolution of the trivial left
    $\Z S\one$-module $\Z$:
    \[
        \Z \xleftarrow{\epsilon} \Z S\one x_{00}
        \xleftarrow{d_1} {\Z S\one}
        \xleftarrow{d_2} \Z S\one x_{00} \oplus \Z S\one x_{01}
        \leftarrow 0
    \]
    Here, $\epsilon(a)=1$ for $a \in S\one x_{00}$,
    $d_1(a)=a(x_{00}-x_{10})$, and
    $d_2(a, b) = a + b = ax_{00} + bx_{01}$,
    so the maps $d_n$
    have the form described above.
    Removing the augmentation and tensoring with $\Z$,
    we find a complex
    \[
        \Z \xleftarrow{0} \Z
        \xleftarrow{\left(\begin{smallmatrix}1 & 1\end{smallmatrix}\right)} \Z \oplus \Z
        \leftarrow 0,
    \]
    which has the homology
    of the 2-sphere $\mathbb{S}^2$. Since $\pi_1(BS) = GS$ is trivial,
    it follows that indeed $BS \simeq \mathbb{S}^2$ by the uniqueness of Moore spaces
    \cite[Example 4.34]{HatcherAT}.
\end{example}

A projective resolution of $\Z$ begins as $\Z \xleftarrow{\epsilon} \Z M e$
for some idempotent $e$, chosen with minimal $|Me|$.
To extend a partial projective resolution
($\Z \leftarrow P_0 \leftarrow \cdots \xleftarrow{d_n} P_{n}$)
by a new map and module ($\xleftarrow{d_{n+1}} P_{n+1}$)
there are two major phases: first, we must compute the kernel of $d_n$,
and second, we must find a map $d_{n+1}$ out of some projective $P_{n+1}$
such that $d_{n+1}(P_{n+1}) = \ker d_n$. The two
phases are described in subsections \ref{sub:kernels}
and \ref{sub:covering} respectively.
Both phases use a newly developed integer linear algebra package
\texttt{mutable\_lattice} \cite{mutable_lattice}, the design of which is described
in subsection \ref{sub:mutable_lattices}.
A further optimization is described in subsection~\ref{sub:splitting_and_caching}.
The final process of tensoring with $\Z$ and computing homology
is described in subsection~\ref{sub:tensoring}.

\subsection{Computing Kernels}\label{sub:kernels}
To compute the kernel of a map of projective $\Z M$-modules
as described above,
we ignore the $\Z M$-module structure and compute a kernel
of a map between free abelian groups.

There is already substantial literature on computing kernels of integer matrices,
many of which have been incorporated
into the SageMath system of \cite{sagemath} for dense integer matrices\footnote{
    See SageMaths's \texttt{Matrix\_integer\_dense.right\_kernel\_matrix()} method.
},
which includes optional delegation
to packages by \cite{PARI2}, \cite{flint}, and \cite{IML}.
However, when compared to alternative packages in the computational ecosystem,
the new \texttt{mutable\_lattice} package proved to be
the fastest at computing kernels
among the alternatives tried when applied
to the relatively sparse matrices with mostly small integer entries
arising in finite monoid homology calculations---see
the performance comparison in Table~\ref{tab:benchmarks}.

The \texttt{mutable\_lattice.relations\_among} method computes
the left kernel of an integer matrix, given as a list of integer row vectors.
Its implementation begins by removing duplicate vectors,
adding relations $(\dots, 1, \dots, -1, \dots)$
to the result that witness any duplication. It then searches
to identify columns in which there is only one nonzero entry
and removes that row and column from consideration.
Next, it partitions the remaining rows and columns
into smaller subproblems based on which rows and columns
intersect along nonzero entries, potentially
reducing the kernel computation to a number of smaller
kernel computations. After these reductions,
the subproblem kernels are computed via standard Gaussian elimination:
an identity matrix
is adjoined on the right to the given list of row vectors,
imagining it as separated by a vertical bar.
Then (generalized) row operations are applied to
put the matrix in Hermite normal form (HNF).
For example, given vectors $(2, 4), (4, 8), (0, 3) \in \Z^2$,
we compute:
\[
    \left(
    \begin{array}{cc|ccc}
        2 & 4 & 1 & 0 & 0 \\
        4 & 8 & 0 & 1 & 0 \\
        0 & 3 & 0 & 0 & 1 \\
    \end{array}
    \right)
    \xrightarrow{\text{row operations}}
    \left(
    \begin{array}{cc|ccc}
        2 & 1 & 1 & 0 & -1 \\
        0 & 3 & 0 & 0 & 1 \\
        \midrule
        0 & 0 & 2 & -1 & 0 \\
    \end{array}
    \right)
\]
We separate the bottom $k$ rows so that the rows
with only zero entries to the left of the vertical bar
are exactly the rows below the separation.
The row operations on the left of the vertical bar
replaced a matrix $A$ by a new matrix $UA$ with the same
row span, where $U$ is an invertible integer matrix representing
the row operations. The matrix $U$ is recorded
to the right of the vertical bar because we applied
the same row operations to the identity matrix.
Because the top left portion of the matrix is in HNF,
its rows are linearly independent, and so $\{w \in \Z^r : w(UA)=0\}$
is spanned by the last $k$ standard basis vectors of $\Z^r$.
Hence, the left kernel of $A$ is
\begin{align*}
    \{w\in \mathbb{Z}^r : wA = 0 \}
    &= \{w'U : w'UA = 0\}
    \\&= \{w' : w'(UA)=0\}\cdot U
    \\&= \langle\text{the last $k$ standard basis vectors of $\Z^r$}\rangle \cdot U
    \\&= \langle\text{the bottom $k$ rows of $U$}\rangle.
\end{align*}
Thus, the rows in the bottom right separation
give a basis for the kernel.
In the displayed example, we found one relation: $2(2,4) -1 (4, 8) + 0(0, 3)=0$.

\subsection{Covering by a New Module}\label{sub:covering}
Once we have computed an ordered $\Z$-basis $\{w_1, \dots, w_k\}$
for the kernel of the outgoing map $d_{n} \colon P_{n} \to P_{n-1}$,
we are left to cover this kernel by the image of some
new map $d_{n+1} \colon P_{n+1} \to P_n$ out of a new projective module $P_{n+1}$.
One inefficient method of doing this is to take $P_{n+1} = (\Z M)^k$
and define $d_{n+1}(x_1, \dots, x_k) \coloneq x_1w_1 + \dots + x_k w_k$.
The image $d_{n+1}$ is then $\Z M w_1 + \cdots +\Z M w_k$, which certainly
contains $\ker d_{n}$, and is a subset of $\ker d_n$
because $\ker d_n$ is closed under left multiplication by $\Z M$
by the left $\Z M$ linearity of $d_n$.

To construct a smaller module $P_{n+1}$,
we can search for a subset
$\{w_{i_1}, \dots, w_{i_r}\} \subseteq \{w_1, \dots, w_k\}$
with the same $\Z M$-span, so then
\begin{equation}
    \label{eq:vecsubset}
    \Z M w_{i_1} + \cdots + \Z M w_{i_r}
    = \ker d_n = \Z w_1 + \cdots + \Z w_{k}.
    \tag{$\star$}
\end{equation}
To search for such a subset, we maintain
a representation of a subgroup
$L \coloneq \Z w_{i_1} + \cdots + \Z w_{i_{r_0}}$ of $\ker d_n$.
For each original basis vector $w_i \in \{w_1, \dots, w_r\}$,
if $w_i$ is not already in $L$ then we add $w_i$
to our solution subset and add the vectors of $\Z M w_i$
to $L$. The subgroup $L$ is stored as a $\Z$-basis
using a \texttt{mutable\_lattice.Lattice} object,
and adding $\Z M w_i$ involves adding each vector $xw_i$
for $x \in M$.
The result of this pass ensures that each
vector in our solution is not in the $\Z M$ span of previously added vectors:
\[
    w_{i_j} \notin \Z M w_{i_1} + \cdots + \Z M w_{i_{j-1}}.
\]
To further ensure that
\[
    w_{i_j} \notin \Z M w_{i_{j+1}} + \cdots + \Z M w_{i_{r}},
\]
we repeat a similar process of filtering
by iterating over the remaining basis in reverse order.
It is possible to further ensure that the chosen subset
is minimal, in the sense that
\[
    w_{i_j} \notin \Z M w_{i_1} + \cdots + \Z M w_{i_{j-1}} + \Z M w_{i_{j+1}} + \cdots + \Z M w_{i_{r}},
\]
but this is more computationally expensive, so we only do this when $k$ is small.
Once we have found an appropriately small subset satisfying (\ref{eq:vecsubset}),
we can further shrink the new module $P_{n+1}$
by finding idempotents $e_1, \dots, e_r \in M$
such that $e_j w_{i_j} = w_{i_j}$, so then
\[
    \Z M e_1 w_{i_1} + \cdots + \Z M e_r w_{i_r}
    = \Z M w_{i_1} + \cdots + \Z M w_{i_r}
    = \ker d_n.
\]
We then use $P_{n+1} = \Z M e_1 \oplus \cdots \oplus \Z M e_r$,
with $d_{n+1}\colon P_{n+1} \to P_n$ defined
as $d_{n+1}(x_1, \dots, x_r) = x_1 w_{i_1} + \cdots + x_r w_{i_r}$,
so $d_{n+1}(P_{n+1}) = \ker d_n$.
Such a choice of idempotents $e_j$ satisfying $e_j w_{i_j} = w_{i_j}$
always exists because we can always take $e_j = 1$,
but we always choose $e_j$ with minimal $|M e_j|$ so the $\Z$-rank
$|Me_1| + \cdots + |M e_r|$ of $P_{n+1}$ is minimized.

The file \texttt{find\_generating\_subset.py} in
the \texttt{fast\_semigroup\_homology} package implements
the process described in this subsection.

\subsection{\texttt{mutable\_lattice}}\label{sub:mutable_lattices}
The covering process described in subsection~\ref{sub:covering}
involves maintaining a representation of a sublattice of some integer lattice $\Z^N$.
An existing implementation\footnote{
    Specifically SageMath's \texttt{FreeModule\_submodule\_pid\_with\_category}
    type prodiced by the expression \texttt{span([vector([1])])}
}
in SageMath was immutable,
so new vectors could not be added without substantial copying.
This led the author to develop the \texttt{mutable\_lattice} package,
available at \cite{mutable_lattice}.

For performance, \texttt{mutable\_lattice} is implemented primarily in the C language,
using Python's C API. Because the integers we use are typically small,
we use a system of tagged pointers to allow efficient storage
of both small integers that fit into one machine word and
pointers to unbounded Python integer objects.

The \texttt{mutable\_lattice} package provides two classes,
\texttt{Vector} and \texttt{Lattice}, as well as the
function \texttt{relations\_among} described in subsection \ref{sub:kernels}.
The \texttt{Vector} class stores a sequence of integers,
and the \texttt{Lattice} class provides a mutable representation
of a subgroup of $\Z^N$. A \texttt{Lattice} object \texttt{L} provides
an efficient \texttt{L.\_\_contains\_\_(v)} method
to check for the presence of a vector in a lattice
and an \texttt{L.add\_vector(v)} method for
replacing a subgroup $L$ by the subgroup $L + \mathbb{Z}v$.

A \texttt{Lattice} object stores a basis in Hermite normal form (HNF).
To add a vector $v$, we first apply row operations between the existing
basis and $v$ to increase the number of leading zeros in $v$,
until the leading nonzero entry of $v$ cannot be reduced---this
nonzero entry becomes a pivot, and the row is inserted into the basis
so that each basis vector has more leading zeros than the previous.
Once inserted, additional row operations are applied so
that the entries above each pivot are normalized to be nonnegative
and less than the pivot value.
A normalization order inspired by \cite{ChouCollins}
is used to maintain smaller integer entries:
we reduce the second-to-last row by the pivot of the last row,
then reduce the third-to-last row by the pivots of the last two
rows, and so on, until all entries above pivots are reduced.
\begin{example}
    To add the last row $(0, 0, 10, 20, 30, 40) \in \Z^6$
    to the lattice spanned by the first four rows below,
    we apply the following process:
    \begin{align*}
        \begin{pmatrix}
            1 & 123 & 1 & 500 & 0 & 3 \\
            0 & 0 & 2 & -10 & 3 & 0 \\
            0 & 0 & 0 &   0 & 5 & 2 \\
            0 & 0 & 0 &   0 & 0 & 4 \\
            \midrule
            0 & 0 & 10 & 20 & 30 & 40 \\
        \end{pmatrix}
        \\\xrightarrow{\text{Operations on new row, then insert}}
        &\begin{pmatrix}
            1 & 123 & 1 & 500 & 0 & 3 \\
            0 & 0 & 2 & -10 & 3 & 0 \\
            0 & 0 & 0 & 70 & 15 & 40 \\
            0 & 0 & 0 &   0 & 5 & 2 \\
            0 & 0 & 0 &   0 & 0 & 4 \\
        \end{pmatrix}
        \\\xrightarrow{\text{Normalize third row by pivots below}}
        &\begin{pmatrix}
            1 & 123 & 1 & 500 & 0 & 3 \\
            0 & 0 & 2 & -10 & 3 & 0 \\
            0 & 0 & 0 &  70 & 0 & 2 \\
            0 & 0 & 0 &   0 & 5 & 2 \\
            0 & 0 & 0 &   0 & 0 & 4 \\
        \end{pmatrix}
        \\\xrightarrow{\text{Normalize second row by pivots below}}
        &\begin{pmatrix}
            1 & 123 & 1 & 500 & 0 & 3 \\
            0 & 0 & 2 &  60 & 3 & 2 \\
            0 & 0 & 0 &  70 & 0 & 2 \\
            0 & 0 & 0 &   0 & 5 & 2 \\
            0 & 0 & 0 &   0 & 0 & 4 \\
        \end{pmatrix}
        \\\xrightarrow{\text{Normalize first row by pivots below}}
        &\begin{pmatrix}
            1 & 123 & 1 & 10 & 0 & 1 \\
            0 & 0 & 2 &  60 & 3 & 2 \\
            0 & 0 & 0 &  70 & 0 & 2 \\
            0 & 0 & 0 &   0 & 5 & 2 \\
            0 & 0 & 0 &   0 & 0 & 4 \\
        \end{pmatrix}
    \end{align*}
\end{example}
We do not use more sophisticated modular
arithmetic algorithms as described in \cite{Storjohann98}
or \cite{spaceHNF}
because the matrices we find in finite monoid homology
calculations are somewhat sparse and have mostly small entries.

To check for containment of a vector $v\in \Z^N$
within a lattice $L\subseteq\Z^N$, we subtract multiples
of the basis vectors of $L$ to increase the number of leading zeros
until either $v$ is the zero vector
(in which case the original $v$ was present in $L$),
or we find a nonzero entry in $v$ that cannot be
made zero with some pivot in the matrix for $L$
(in which case the original $v$ was absent from $L$).
A fast initial check can also be performed to check whether
$v$ has any nonzero entries at indexes where the basis of $L$
has only zero entries, in which case $v$ is certainly not present.

\subsection{Splitting and Caching}\label{sub:splitting_and_caching}
When constructing projective resolutions of the
form described so far, it is sometimes the case
that the kernel of a map
\[
    d \colon \bigoplus_{j=1}^m \Z M e_j \to \bigoplus_{i=1}^{m'} \Z M f_i
\]
splits along the given summands: we can partition
$\{1, \dots, m\} = A_1 \sqcup \cdots \sqcup A_q$
so that $\ker d = \iota_{A_1}(K_1) + \cdots + \iota_{A_q}(K_q)$,
where each $\iota_{A_\l}$ is the inclusion
\[
    \iota_{A_\l} \colon \bigoplus_{j \in A_\l} \Z M e_j \hookrightarrow \bigoplus_{j=1}^{m} \Z M e_j,
\]
and each $K_\l$ is some submodule of $\bigoplus_{j \in A_\l} \Z M e_j$.
When such splitting occurs,
we can solve the the corresponding covering
problem for each $K_\l$ individually, since each $\iota_{A_1}(K_\l)$ is closed
under the $\Z M$ action. The \texttt{decompose} method
on \texttt{Lattice} objects is used to find such a splitting.
\begin{example}\label{ex:rect32}
    Consider the $3$-by-$2$ rectangular band
    $S = \{x_{ij} : i \in \{0,1,2\}, j \in \{0,1\}\}\}$,
    which has operation defined by $x_{ij}x_{i'j'} = x_{ij'}$.
    We build the following projective resolution
    of the trivial left $\Z S\one$-module $\Z$:
    \[
        \Z \xleftarrow{\epsilon}
        \Z S\one x_{00}
        \xleftarrow{d_1}
        \left[\begin{matrix}
        \Z S\one & \xleftarrow{d_{2,0}} & \Z S\one x_{00} \oplus \Z S\one x_{01} & \leftarrow & 0
        \\
        \oplus & \phantom{\xleftarrow{d_{2,0}}} & \oplus & {} & \oplus
        \\
        \Z S\one & \xleftarrow{d_{2,1}} & \Z S\one x_{00} \oplus \Z S\one x_{01} & \leftarrow & 0
        \end{matrix}\right.
    \]
    The maps here are defined as $\epsilon(a) = 1$
    for $a \in S\one x_{00}$,
    $d_1(a, b) = a(x_{00}-x_{20}) + b(x_{10}-x_{20})$,
    and $d_{2,i}(a,b)= a + b = ax_{00} + b x_{01}$.
    In this example, $\ker d_1$ splits along the two
    copies of $\Z S\one$, so the two summands
    of the kernel can be covered separately
    by $d_{2,0}$ and $d_{2,1}$.
    After removing the augmentation and tensoring with $\Z$,
    we find the complex
    \[
        \Z
        \xleftarrow{0}
        \left[\begin{matrix}
        \Z & \xleftarrow{\left(\begin{smallmatrix}1 & 1\end{smallmatrix}\right)} & \Z^2 & \leftarrow & 0
        \\
        \oplus & \phantom{\xleftarrow{\left(\begin{smallmatrix}1 & 1\end{smallmatrix}\right)}} & \oplus & {} & \oplus
        \\
        \Z & \xleftarrow{\left(\begin{smallmatrix}1 & 1\end{smallmatrix}\right)} & \Z^2 & \leftarrow & 0
        \end{matrix}\right.
    \]
    which has the homology $H_2(BS) = \Z^2$ and $\tilde H_i(BS)=0$ for $i \neq 2$.
    Because $\pi_1(BS)=GS=1$, the space $BS$ is a Moore space of type $M(\Z^2, 2)$,
    so in fact $BS \simeq \mathbb{S}^2 \vee \mathbb{S}^2$ by the
    uniqueness of Moore spaces.
\end{example}

The preceding example may be optimized further by noting that
the two summands of $\ker d_1$ are the same subset $\Z S \subset \Z S\one$,
and so we can re-use previously created parts of the resolution
by  maintaining a cache of which subsets of which projective modules
are covered by which new maps.
The \texttt{fast\_semigroup\_homology} package applies this optimization,
which is useful in the following notable case:
\begin{example}
    For the cyclic group $C_k=\{1, a, \cdots,a^{k-1}\}$,
    we have the following classical projective resolution
    of the trivial left $\Z C_k$-module $\Z$:
    \[
        \Z \xleftarrow{\epsilon} \Z C_k
        \xleftarrow{\cdot(a - 1)} \Z C_k
        \xleftarrow{\cdot N} \Z C_k
        \xleftarrow{\cdot(a - 1)} \Z C_k
        \xleftarrow{\cdot N} \Z C_k
        \xleftarrow{\cdot(a - 1)} \Z C_k
        \leftarrow \cdots
    \]
    Here $N = 1+a + \cdots + a^{k-1}$.
    However, when constructing this resolution,
    the \texttt{fast\_semigroup\_homology} package
    identifies that $\ker(\cdot N)$ and $\ker \epsilon$
    are the same subset of $\Z C_k$, so we can re-use the same
    module in dimension $1$ as the module in dimension $3$.
    The result is that the above infinite resolution can
    be represented as the following finite graph with a cycle:
    \[
    \begin{tikzpicture}
        [every node/.style={scale=1.0,inner sep=0.05cm,outer sep=0}]
        \node (start)
            {$\Z \xleftarrow{\epsilon} \Z C_k$};
        \node (rest)
            [right=1cm of start]
            {$\Z C_k \xleftarrow{N} \Z C_k$};
        \draw[->] ($(rest.west)+(0,-0.081)$) -- ($(start.east)+(0,-0.081)$)
            node[midway, above=0.0cm] {\scriptsize$a-1$};
        \draw[->, rounded corners] ($(rest.west)+(0,-0.081)$)
                -- +(-0.4, 0)
                -- ($(rest.south west)+(-0.4, -0.2)$)
                -- ($(rest.south east)+(+0.4, -0.2)$)
                -- ($(rest.east)+(+0.4, -0.081)$)
                -- ($(rest.east)+(0,-0.081)$)
                ;
    \end{tikzpicture}
    \]
    Removing the augmentation and tensoring with $\Z$ results in the complex
    \[
    \begin{tikzpicture}
        [every node/.style={scale=1.0,inner sep=0.05cm,outer sep=0}]
        \node (start)
            {$\Z$};
        \node (rest)
            [right=1cm of start]
            {$\Z \xleftarrow{k} \Z$};
        \draw[->] ($(rest.west)+(0,-0.081)$) -- ($(start.east)+(0,-0.081)$)
            node[midway, above=0.05cm] {\scriptsize{$0$}};
        \draw[->, rounded corners] ($(rest.west)+(0,-0.081)$)
                -- +(-0.4, 0)
                -- ($(rest.south west)+(-0.4, -0.2)$)
                -- ($(rest.south east)+(+0.4, -0.2)$)
                -- ($(rest.east)+(+0.4, -0.081)$)
                -- ($(rest.east)+(0,-0.081)$)
                ;
    \end{tikzpicture}
    \]
    from which we compute the usual homology groups
    $H_0(C_k)=\Z$, $H_i(C_k)=C_k$ for odd $i$,
    and $H_i(C_k)=0$ for even $i \geq 2$.
\end{example}
The implementation of this technique of caching and re-using
parts of a resolution is implemented by
\texttt{projective\_resolution.py}
in the \texttt{fast\_semigroup\_homology} package.
More examples of this technique
are presented in Section~\ref{sec:results}, particularly Example~\ref{ex:2pow},
in which the technique provided an exponential speedup.

\subsection{Tensoring and Homology}\label{sub:tensoring}

Once we have found a sufficiently long partial projective resolution $P_*$
of the trivial left $\Z M$-module $\Z$,
to compute the homology groups $H_i(M) = \Tor_i^{\Z M}(\Z, \Z)$,
we need to first find the groups and boundary maps
in the chain complex $\Z \tensor_{\Z M} P_*$.

For the projective modules $P_n = \bigoplus_{j=1}^m \Z M e_j$ we use,
because $\Z \tensor_{\Z M} \Z M e \cong \Z$
via the isomorphism $1 \otimes x \mapsto 1$ for $x \in Me$,
we have $\Z \tensor_{\Z M} P_n \cong \Z^m$.
Working along the direct sums,
each map
$d \colon \bigoplus_{j=1}^m \Z M e_j \to \bigoplus_{i=1}^{m'} \Z M f_j$
defined by $p_i(d(\iota_j(x))) = xa_j^i$,
descends to the tensored map
$\id \tensor d = \id_\Z \tensor_{\Z M} d$
which may be defined by
$p_i((\id \tensor d)(\iota_j(1 \tensor x)) = (1 \tensor x)\epsilon(a_i^j)$,
where $\epsilon \colon \Z M \to \Z$ has $\epsilon(x)=1$ for $x \in M$.
Thus, the matrix defining the tensored map $\Z^{m} \to \Z^{m'}$
representing $\id \tensor d$ has entries $\epsilon(a_i^j)$.
The boundary matrices in the tensored complex $\Z \tensor_{\Z M} P_*$
can therefore be computed by summing the coefficients
within each $\Z M$ element entry of the boundary matrices for $P_*$.
If there is a splitting
in the resolution as described in subsection~\ref{sub:splitting_and_caching},
we still have the same splitting after tensoring with $\Z$,
as seen in Example~\ref{ex:rect32}.

To finally compute the homology of a chain complex,
the standard technique uses Smith normal forms \cite[p.60]{rotmanAT}:
given a chain complex
\[
    \Z^{m''} \xleftarrow{d'} \Z^{m'} \xleftarrow{d} \Z^m,
\]
the homology at the middle term is given by
\[
    \ker d' / \im d \cong \Z^{m' - \rank{d'} - \rank{d}} \oplus \Z/r_1 \oplus \cdots \oplus \Z/r_k,
\]
where $r_1, \dots, r_k$ are the nonzero diagonal entries of
some (not necessarily square) diagonal matrix representing $d$
in some pair of bases. For a complex
\[
    \Z^{m''}
    \xleftarrow{d'}
    \bigoplus_{i=1}^q \Z^{m'_i}
    \xleftarrow{d=\bigoplus_{i=1}^q d_i}
    \bigoplus_{i=1}^q \Z^{m_i},
\]
with splitting
structure as described in subsection~\ref{sub:splitting_and_caching},
we can find such diagonal entries of each summand map $d_i$ independently,
then concatenate the lists to find the diagonal entries of the sum map $d$.
This allows re-using the invariant factors $r_1, \dots, r_k$
for the tensored outgoing map from each module after computing them once,
even if the module is re-used
several times in the projective resolution.
The diagonal matrices are traditionally normalized
to satisfy the Smith normal form divisibility constraint $r_1 \mid \dots \mid r_k$,
but this can be enforced after concatenation by repeatedly
replacing $\Z/a \oplus \Z/b$
with the isomorphic abelian group $\Z/\gcd(a,b) \oplus \Z/\lcm(a,b)$.

The code used to find the invariant factors $r_1, \dots, r_k$
of a matrix is implemented by the \texttt{nonzero\_invariants}
method on \texttt{Lattice} objects in the \texttt{mutable\_lattice} package,
which is based on the technique of alternating transposition with HNF
calculations described in \cite{HNFSNF}.
The code used to convert lists of invariants into homology
groups is implemented in the files \texttt{projective\_resolution.py}
and \texttt{normalized\_invariants.py} in the \texttt{fast\_semigroup\_homology}
package.

\section{Finding Semigroups and Monoids}\label{sec:finding}
The Smallsemi GAP package \cite[]{smallsemi}
already provides multiplication tables for all semigroups of order at most 8.
However, the numbers of semigroups of order $n$ grows extremely quickly as $n$ grows
(see Table~\ref{tab:counts_table}),
and databases of semigroups of larger orders have not been publicly available.
Fortunately, the vast majority of known semigroups have
some idempotent $e$ satisfying the hypothesis of Theorem~\ref{thm:eSeEquivalent},
so there is a monoid $eSe$ with fewer elements having a homotopy equivalent
classifying space.
In searching for finite semigroups with unique classifying spaces,
we can therefore limit our search to the much smaller set of
only those semigroups of a particular order with no non-identity idempotent $e$
satisfying $eSe=eS$ or $eSe=Se$. If this was the only restriction,
our search would often include both a semigroup $S$ of order $n$
and the monoid $S\one$ of order $n+1$,
which have homotopy equivalent classifying spaces; instead
we restrict our search to finite monoids, so each non-monoid semigroup's
classifying space is found as the classifying space of the monoid $S\one$.
Finally, for each finite semigroup $S$, we examine
only one of $S$ or $S^{\op}$, since these have
homotopy equivalent classifying spaces as well.

We follow techniques described in \cite{DistlerThesis}
and in \cite{mon8910} to generate finite semigroups/monoids
using the \texttt{minion} solver \cite[]{minion_paper,minion}.
The program used to generate semigroups/monoids is available
in the \texttt{semisearch} repository of \cite[]{semisearch}.
The numbers of semigroups and monoids
found are are summarized in Table~\ref{tab:counts_table}.

Because of Corollary~\ref{cor:kThin},
we generate finite semigroups or monoids
based on the possible structures of the minimal ideal $K(S)$.
For example, the smallest $K(S)$ for which the
hypothesis of Corollary~\ref{cor:kThin} does not hold is a $2$-by-$2$
rectangular band $\{x_{00}, x_{01}, x_{10}, x_{11}\}$
with operation $x_{ij}x_{i'j'} = x_{ij'}$. By fixing this minimal
ideal as the first four elements of a multiplication table,
it suffices to search only for multiplication tables in the format of
Table~\ref{tab:KS_is_rect22}.
In general,
we can enforce the fact that $K(S)$ is the minimal ideal
by choosing $K(S)$ to be simple and applying the constraint
that it absorbs multiplication by other elements.
We can also enforce prior knowledge of the minimal
left and right ideals of $K(S)$ based on the
Rees matrix semigroup structure of $K(S)$.

We also follow Section 3.1 of \cite{DistlerThesis} in separating
the search space by the isomorphism type of the diagonal
of the multiplication table, and we follow \cite{mon8910}
in also separating by the group of units when searching for monoids.
Separating the search space into many parts with more structure has two main purposes:
(1) it allows the individual search processes to be easily run in parallel,
and (2) it reduces the number of permutations of the elements that fix that the structure,
so there are fewer isomorphic multiplication tables to compare against to remove duplicates.
If there are only a small number of structure-stabilizing permutations, we include
symmetry-breaking constraints directly into the \texttt{minion} problem
by requiring solution multiplication tables to be lexicographically minimal
among isomorphic tables with the given structure.
For larger stabilizers, we must apply a post-filtering step to individually
compare against each such isomorphic table. This follows
the ``Complete Symmetry-Breaking'' versus ``Isomorph Rejection''
distinction defined in \cite{mon8910}.

\section{Computational Results}\label{sec:results}

The \texttt{fast\_semigroup\_homology} software from Section~\ref{sec:fast}
was used to calculate the homology of finite semigroups and monoids generated by the
by the \texttt{semisearch} program described in Section~\ref{sec:finding}.
The results of these calculations are available
in the ``results'' folder of \cite{fast_semigroup_homology},
and include homology for all semigroups of order at most 8,
along with homology for all monoids not satisfying the
hypothesis of Theorem~\ref{thm:eSeEquivalent} of order at most 11,
and a further restricted class of monoids of order at most 16.
The entire folder took several hours to generate,
and specific runtimes are given in the \texttt{.md} files of results.
Limited subsets of the results are presented
in Table~\ref{tab:semigroup_homology_tables}
and Table~\ref{tab:monoid_homology_tables}.

The remainder of this section highlights a selection of finite semigroups
found to have interesting properties,
as described by Table~\ref{tab:false_conjectures}.
Each projective resolution calculation below is verified
by a test case in the \texttt{test\_projective\_resolution.py} file
of the \texttt{fast\_semigroup\_homology} package.

\begin{example}\label{ex:no_zero_contractible}
    Consider the semigroup $S=\{x_{00}, x_{01}, x_{10}, x_{11}, q\}$
    with operation defined by $x_{ij}x_{i'j'} = x_{ij'}$
    and $x_{ij}q=x_{i0}$
    and $qx_{ij}=x_{ij}$ and $qq=q$. This can be written as
    the following multiplication table:
    \[
    \text{
    \footnotesize
    \tabcolsep=0.1cm
        \begin{tabular}{c|ccccc}
        $\cdot$  & 0 & 1 & 2 & 3 & 4 \\ \hline
        0        & 0 & 1 & 0 & 1 & 0 \\
        1        & 0 & 1 & 0 & 1 & 0 \\
        2        & 2 & 3 & 2 & 3 & 2 \\
        3        & 2 & 3 & 2 & 3 & 2 \\
        4        & 0 & 1 & 2 & 3 & 4 \\
        \end{tabular}
    }
    \]
    This semigroup $S$ has $BS$ contractible but no left or
    right zero element.
    \begin{proof}
        Because $Sq=\{x_{00}, x_{10}, q\}=qSq$, Theorem~\ref{thm:eSeEquivalent}
        implies that $BS \simeq B(qSq)$, and the monoid $qSq$
        has a left zero element $x_{00}$, so $B(qSq)$ is contractible.
    \end{proof}
    There are also semigroups with 7 elements (or monoids with 8 elements)
    that have trivial homology but do not satisfy
    the hypothesis of Theorem~\ref{thm:eSeEquivalent}---see
    the folder of complete results.
\end{example}

\begin{example}\label{ex:0ZZZ}
    Consider the semigroup $S=\{x_{00}, x_{01}, x_{10}, x_{11}, y_{00}\}$
    with operation defined by $x_{ij}x_{i'j'} = x_{ij'}$
    and $y_{00}x_{ij} = x_{0j}$
    and $x_{ij}y_{00}=x_{i0}$
    and $y_{00}y_{00} = x_{00}$. This can be written as the following
    multiplication table:
    \[
    \text{
    \footnotesize
    \tabcolsep=0.1cm
        \begin{tabular}{c|ccccc}
        $\cdot$  & 0 & 1 & 2 & 3 & 4 \\ \hline
        0        & 0 & 1 & 0 & 1 & 0 \\
        1        & 0 & 1 & 0 & 1 & 0 \\
        2        & 2 & 3 & 2 & 3 & 2 \\
        3        & 2 & 3 & 2 & 3 & 2 \\
        4        & 0 & 1 & 0 & 1 & 0 \\
        \end{tabular}
    }
    \]
    This semigroup $S$ has homology $H_1(S)=0$ but
    $H_i(S)=\Z$ for all $i\geq 2$.
    \begin{proof}
        The \texttt{fast\_semigroup\_homology} package
        finds the following repeating
        projective resolution of the trivial $\Z S\one$-module $\Z$:
        \[
        \begin{tikzpicture}
            [every node/.style={scale=1.0,inner sep=0.02cm,outer sep=0}]
            \node (aug) {$\Z$};
            \node (dimzero)
                [right=0.75cm of aug]
                {$\Z S\one x_{00}$};
            \node (dimone)
                [right=1.05cm of dimzero]
                {$\Z S\one$};
            \node (dimtwoOplus)
                [right=2cm of dimone]
                {$\oplus$};
            \node (dimtwoTop)
                [above=0cm of dimtwoOplus]
                {$\Z S\one x_{00}$};
            \node (dimtwoBottom)
                [below=0cm of dimtwoOplus]
                {$\Z S\one$};
            \node (dimthreeOplus)
                [right=2cm of dimtwoBottom]
                {$\oplus$};
            \node (dimthreeTop)
                [above=0cm of dimthreeOplus]
                {$\Z S\one x_{10}$};
            \node (dimthreeBottom)
                [below=0cm of dimthreeOplus]
                {$\Z S\one$};
            \coordinate (dimOneBottomLeftCorner) at ($(dimone.south west)+(-0.03,-0.04)$);
            \coordinate (dimOneTopLeftCorner) at (dimone.north west-|dimOneBottomLeftCorner);
            \draw[thick] ($(dimOneBottomLeftCorner)+(0.2,0)$)
                -- (dimOneBottomLeftCorner)
                -- (dimOneTopLeftCorner)
                -- +(0.2, 0);
            \node (dimOneLeft)
                [circle, fill, minimum size=0.1cm]
                at ($(dimOneBottomLeftCorner)!0.5!(dimOneTopLeftCorner)$)
                {};
            \coordinate (dimTwoTopLeftCorner) at ($(dimtwoTop.north west)+(-0.03,0)$);
            \coordinate (dimTwoBottomLeftCorner) at (dimtwoBottom.south west-|dimTwoTopLeftCorner);
            \draw[thick] ($(dimTwoBottomLeftCorner)+(0.2,0)$)
                -- (dimTwoBottomLeftCorner)
                -- (dimTwoTopLeftCorner)
                -- +(0.2, 0);
            \node (dimTwoLeft)
                [circle, fill, minimum size=0.1cm]
                at (dimOneLeft-|dimTwoBottomLeftCorner)
                {};
            \coordinate (dimThreeTopLeftCorner) at ($(dimthreeTop.north west)+(-0.03,0)$);
            \coordinate (dimThreeBottomLeftCorner) at (dimthreeBottom.south west-|dimThreeTopLeftCorner);
            \draw[thick] ($(dimThreeBottomLeftCorner)+(0.2,0)$)
                -- (dimThreeBottomLeftCorner)
                -- (dimThreeTopLeftCorner)
                -- +(0.2, 0);
            \node (dimThreeLeft)
                [circle, fill, minimum size=0.1cm]
                at ($(dimThreeBottomLeftCorner)!0.5!(dimThreeTopLeftCorner)$)
                {};
            \draw[->] (dimzero.west) -- (aug.east)
                node[midway, above=0.1cm] {\scriptsize$\epsilon$};
            \draw[->] (dimOneLeft) -- (dimOneLeft-|dimzero.east)
                node[midway, above=0.1cm] {\scriptsize$\partial_1$};
            \draw[->] (dimTwoLeft) -- (dimTwoLeft-|dimone.east)
                node[midway, above=0.1cm] {\scriptsize$\partial_2$};
            \draw[->] (dimThreeLeft) -- (dimThreeLeft-|dimtwoBottom.east)
                node[midway, above=0.1cm] {\scriptsize$\partial_3$};
            \draw[->, rounded corners] (dimTwoLeft)
                -- +(-0.6, -0.0)
                -- ($(dimTwoBottomLeftCorner)+(-0.6, -0.5)$)
                -- ($(dimthreeBottom.south east)+(+0.4, -0.2)$)
                -- ($(dimthreeBottom.east)+(+0.4, 0)$)
                -- ($(dimthreeBottom.east)$)
                ;
        \end{tikzpicture}
        \]
        Above, $\partial_1(a)=a(x_{00}-x_{10})$, $\partial_2(a, b) = ax_{01} + by_{00}$,
        and $\partial_3(a, b) = a(x_{10}-x_{11}) + b(x_{00} - y_{00})$,
        and the arrows indicate that $\ker \partial_2 = 0 \oplus \im \partial_3$
        and $\ker \partial_3 = 0 \oplus \im \partial_2$.
        Removing the augmentation and tensoring with $\Z$ gives the following
        complex, which has the stated homology:
        \[
        \begin{tikzpicture}
            [every node/.style={scale=1.0,inner sep=0.02cm,outer sep=0}]
            \node (dimzero)
                {$\Z$};
            \node (dimone)
                [right=1.05cm of dimzero]
                {$\Z$};
            \node (dimtwoOplus)
                [right=2cm of dimone]
                {$\oplus$};
            \node (dimtwoTop)
                [above=0cm of dimtwoOplus]
                {$\Z$};
            \node (dimtwoBottom)
                [below=0cm of dimtwoOplus]
                {$\Z$};
            \node (dimthreeOplus)
                [right=2cm of dimtwoBottom]
                {$\oplus$};
            \node (dimthreeTop)
                [above=0cm of dimthreeOplus]
                {$\Z$};
            \node (dimthreeBottom)
                [below=0cm of dimthreeOplus]
                {$\Z$};
            \coordinate (dimOneBottomLeftCorner) at ($(dimone.south west)+(-0.03,-0.04)$);
            \coordinate (dimOneTopLeftCorner) at ($(dimone.north west)+(-0.03,+0.04)$);
            \draw[thick] ($(dimOneBottomLeftCorner)+(0.2,0)$)
                -- (dimOneBottomLeftCorner)
                -- (dimOneTopLeftCorner)
                -- +(0.2, 0);
            \node (dimOneLeft)
                [circle, fill, minimum size=0.1cm]
                at ($(dimOneBottomLeftCorner)!0.5!(dimOneTopLeftCorner)$)
                {};
            \coordinate (dimTwoTopLeftCorner) at ($(dimtwoTop.north west)+(-0.03,+0.04)$);
            \coordinate (dimTwoBottomLeftCorner) at ($(dimtwoBottom.south west)+(-0.03,-0.04)$);
            \draw[thick] ($(dimTwoBottomLeftCorner)+(0.2,0)$)
                -- (dimTwoBottomLeftCorner)
                -- (dimTwoTopLeftCorner)
                -- +(0.2, 0);
            \node (dimTwoLeft)
                [circle, fill, minimum size=0.1cm]
                at (dimOneLeft-|dimTwoBottomLeftCorner)
                {};
            \coordinate (dimThreeTopLeftCorner) at ($(dimthreeTop.north west)+(-0.03,+0.04)$);
            \coordinate (dimThreeBottomLeftCorner) at ($(dimthreeBottom.south west)+(-0.03,-0.04)$);
            \draw[thick] ($(dimThreeBottomLeftCorner)+(0.2,0)$)
                -- (dimThreeBottomLeftCorner)
                -- (dimThreeTopLeftCorner)
                -- +(0.2, 0);
            \node (dimThreeLeft)
                [circle, fill, minimum size=0.1cm]
                at ($(dimThreeBottomLeftCorner)!0.5!(dimThreeTopLeftCorner)$)
                {};
            \draw[->] (dimOneLeft) -- (dimOneLeft-|dimzero.east)
                node[midway, above=0.1cm] {\scriptsize$0$};
            \draw[->] (dimTwoLeft) -- (dimTwoLeft-|dimone.east)
                node[midway, above=0.1cm] {\scriptsize$\begin{pmatrix}1 & 1\end{pmatrix}$};
            \draw[->] (dimThreeLeft) -- (dimThreeLeft-|dimtwoBottom.east)
                node[midway, above=0.1cm] {\scriptsize$0$};
            \draw[->, rounded corners] (dimTwoLeft)
                -- +(-0.6, -0.0)
                -- ($(dimTwoBottomLeftCorner)+(-0.6, -0.4)$)
                -- ($(dimthreeBottom.south east)+(+0.4, -0.2)$)
                -- ($(dimthreeBottom.east)+(+0.4, 0)$)
                -- ($(dimthreeBottom.east)$)
                ;
        \end{tikzpicture}
        \qedhere
        \]
    \end{proof}
\end{example}

\begin{example}\label{ex:2pow}
    Extending the previous example,
    consider the semigroup $S = \{x_{00}, x_{01}, x_{10}, x_{11}, y_{00}, z_{00}\}$
    with operation defined by $a_{ij}a'_{i'j'} = x_{ij'}$
    for any $a,a'\in\{x,y,z\}$. This can be written as the following multiplication
    table:
    \[
    \text{
    \footnotesize
    \tabcolsep=0.1cm
        \begin{tabular}{c|cccccc}
        $\cdot$  & 0 & 1 & 2 & 3 & 4 & 5 \\ \hline
        0        & 0 & 1 & 0 & 1 & 0 & 0 \\
        1        & 0 & 1 & 0 & 1 & 0 & 0 \\
        2        & 2 & 3 & 2 & 3 & 2 & 2 \\
        3        & 2 & 3 & 2 & 3 & 2 & 2 \\
        4        & 0 & 1 & 0 & 1 & 0 & 0 \\
        5        & 0 & 1 & 0 & 1 & 0 & 0 \\
        \end{tabular}
    }
    \]
    This semigroup $S$ has $H_1(S)=0$ and $H_i(S)=\Z^{2^{i-2}}$
    for $i \geq 2$.
    \begin{proof}
    The \texttt{fast\_semigroup\_homology} package
    finds a repeating projective resolution for $\Z$
    as a $\Z S\one$-module, of which the following is a simplified version:
    \[
        \begin{tikzpicture}
            [every node/.style={scale=0.8,inner sep=0.02cm,outer sep=0}]
            \node (aug) {$\Z$};
            \node (dimzero)
                [right=0.5cm of aug]
                {$\Z S\one x_{00}$};
            \node (dimone)
                [right=0.7cm of dimzero]
                {$\Z S\one$};
            \node (dimtwoOplus)
                [right=2.2cm of dimone]
                {$\oplus$};
            \node (dimtwoTop)
                [above=0cm of dimtwoOplus]
                {$\Z S\one x_{00}$};
            \node (dimtwoBottom)
                [below=0cm of dimtwoOplus]
                {$\Z S\one x_{00} \oplus \Z S\one \oplus \Z S\one$};
            \node (dimthreeOplusMiddle)
                [right=2cm of dimtwoBottom]
                {$\oplus$};
            \node (dimthreeTopMiddle)
                [above=0cm of dimthreeOplusMiddle]
                {$\Z S\one x_{00}$};
            \node (dimthreeOplusTop)
                [above=0cm of dimthreeTopMiddle]
                {$\oplus$};
            \node (dimthreeTop)
                [above=0cm of dimthreeOplusTop]
                {$\Z S\one x_{00}$};
            \node (dimthreeBottomMiddle)
                [below=0cm of dimthreeOplusMiddle]
                {$\Z S\one x_{00} \oplus \Z S\one \oplus \Z S\one$};
            \node (dimthreeOplusBottom)
                [below=0cm of dimthreeBottomMiddle]
                {$\oplus$};
            \node (dimthreeBottom)
                [below=0cm of dimthreeOplusBottom]
                {$\Z S\one x_{00} \oplus \Z S\one \oplus \Z S\one$};
            \coordinate (dimOneBottomLeftCorner) at ($(dimone.south west)+(-0.03,-0.04)$);
            \coordinate (dimOneTopLeftCorner) at (dimone.north west-|dimOneBottomLeftCorner);
            \draw[thick] ($(dimOneBottomLeftCorner)+(0.2,0)$)
                -- (dimOneBottomLeftCorner)
                -- (dimOneTopLeftCorner)
                -- +(0.2, 0);
            \node (dimOneLeft)
                [circle, fill, minimum size=0.1cm]
                at ($(dimOneBottomLeftCorner)!0.5!(dimOneTopLeftCorner)$)
                {};
            \coordinate (dimTwoBottomLeftCorner) at ($(dimtwoBottom.south west)+(-0.03,0)$);
            \coordinate (dimTwoTopLeftCorner) at (dimtwoTop.north west-|dimTwoBottomLeftCorner);
            \draw[thick] ($(dimTwoBottomLeftCorner)+(0.2,0)$)
                -- (dimTwoBottomLeftCorner)
                -- (dimTwoTopLeftCorner)
                -- +(0.2, 0);
            \node (dimTwoLeft)
                [circle, fill, minimum size=0.1cm]
                % at ($(dimTwoBottomLeftCorner)!0.5!(dimTwoTopLeftCorner)$)
                at (dimOneLeft-|dimTwoBottomLeftCorner)
                {};
            \coordinate (dimThreeBottomLeftCorner) at (dimthreeBottom.south west);
            \coordinate (dimThreeTopLeftCorner) at (dimthreeTop.north west-|dimThreeBottomLeftCorner);
            \draw[thick] ($(dimThreeBottomLeftCorner)+(0.2,0)$)
                -- (dimThreeBottomLeftCorner)
                -- (dimThreeTopLeftCorner)
                -- +(0.2, 0);
            \node (dimThreeLeft)
                [circle, fill, minimum size=0.1cm]
                at ($(dimThreeBottomLeftCorner)!0.5!(dimThreeTopLeftCorner)$)
                {};
            \draw[->] (dimzero.west) -- (aug.east)
                node[midway, above=0.1cm] {\scriptsize$\epsilon$};
            \draw[->] (dimOneLeft) -- (dimOneLeft-|dimzero.east)
                node[midway, above=0.1cm] {\scriptsize$\partial_1$};
            \draw[->] (dimTwoLeft) -- (dimTwoLeft-|dimone.east)
                node[midway, above=0.1cm] {\scriptsize$\partial_2$};
            \draw[->] (dimThreeLeft) -- (dimThreeLeft-|dimtwoBottom.east)
                node[midway, above=0.1cm] {\scriptsize$\partial_3$};
            \draw[->, rounded corners] (dimThreeLeft)
                -- +(-0.2, -0.0)
                -- ($(dimthreeBottom.south west)+(-0.2, -0.2)$)
                -- ($(dimthreeBottom.south east)+(+0.4, -0.2)$)
                -- ($(dimthreeBottom.east)+(+0.4, 0)$)
                -- ($(dimthreeBottom.east)$)
                ;
            \draw[->, rounded corners] (dimThreeLeft)
                -- +(-0.4, -0.0)
                -- ($(dimthreeBottom.south west)+(-0.4, -0.4)$)
                -- ($(dimthreeBottom.south east)+(+0.6, -0.4)$)
                -- ($(dimthreeBottomMiddle.east)+(+0.6, 0)$)
                -- ($(dimthreeBottomMiddle.east)$)
                ;
        \end{tikzpicture}
    \]
    Above, $\partial_1(a) = a(x_{00}-x_{10})$,
    $\partial_2(a,b,c,d) = ax_{01} + bx_{00} + cy_{00} + dz_{00}$,
    and \[
        \partial_3(a,b,c,d,e,f,g,h)
        =\begin{pmatrix}
            (a + b + c + d + e + f + g + h)x_{00} \\
            -ax_{01} - cx_{00} - d y_{00} - ez_{00} \\
            -b x_{01} - fx_{00} - g y_{00} - h z_{00}
        \end{pmatrix}.
    \]
    The arrows indicate that $\ker \partial_2 = 0 \oplus \im \partial_3$
    and $\ker \partial_3 = 0 \oplus 0 \oplus \im \partial_3 \oplus \im \partial_3$.
    The tensored boundary maps are as follows: $\partial_1 \times_{\Z S\one} \id_{\Z}$ is the zero map,
    $\partial_2 \times_{\Z S\one} \id_{\Z}$ is summation $\Z^4 \to \Z$,
    and $\partial_3 \times_{\Z S\one} \id_{\Z}$ is represented by the matrix
    \[
        \left(
        \begin{smallmatrix}
            1 & 1 & 1 & 1 & 1 & 1 & 1 &  1 \\
            -1 & 0 & -1 & -1 & -1 & 0 & 0 & 0 \\
            0 & -1 & 0 & 0 & 0 & -1 & -1 & -1
        \end{smallmatrix}
        \right),
    \]
    which has $[1,1]$ as its nonzero Smith normal form invariant factors.
    The result is that $H_1(S)=0$, $H_2(S)=\Z$,
    and by expanding a tree of $2^{i-2}$ copies of the last bracketed module
    in each dimension $i\geq 2$,
    we see $H_i(S) = (H_2(S))^{2^{i-2}} = \Z^{2^{i-2}}$.
    \end{proof}
\end{example}
Note that the exponential growth of the homology
(and therefore cohomology by the universal coefficient theorem)
above implies that the cohomology ring cannot have a finite set
of generators: the set of all monomials in a finite set of generators grows only
polynomially as the degree increases, so the dimension of the span of these monomials
would only grow polynomially as well.

\begin{example}\label{ex:Moore_C2_3}
    Consider the semigroup $S = \{x_0, \dots, x_8\}$
    with operation defined by $x_ix_j=x_k$ for every entry $k$ in row $i$ and column $j$
    in the following table\footnote{
        This is the table at index \texttt{57244}
        of the \texttt{order10.hdf5} file in the
        ``\texttt{monoids no monoid 1sided ideals by min ideal and diagonal and units}''
        results folder from the \texttt{semisearch} repository,
        with the adjoined identity removed.
    }:
    \[
    \text{
    \footnotesize
    \tabcolsep=0.1cm
        \begin{tabular}{c|ccccccccc}
        $\cdot$  & 0 & 1 & 2 & 3 & 4 & 5 & 6 & 7 & 8 \\ \hline
        0        & 0 & 1 & 0 & 1 & 0 & 0 & 0 & 0 & 0 \\
        1        & 0 & 1 & 0 & 1 & 0 & 0 & 0 & 1 & 1 \\
        2        & 2 & 3 & 2 & 3 & 2 & 2 & 2 & 2 & 2 \\
        3        & 2 & 3 & 2 & 3 & 2 & 2 & 2 & 3 & 3 \\
        4        & 0 & 1 & 0 & 1 & 0 & 4 & 4 & 0 & 0 \\
        5        & 0 & 1 & 2 & 3 & 0 & 5 & 5 & 0 & 0 \\
        6        & 0 & 1 & 2 & 3 & 0 & 6 & 6 & 0 & 0 \\
        7        & 0 & 1 & 0 & 1 & 4 & 0 & 4 & 7 & 8 \\
        8        & 0 & 1 & 0 & 1 & 4 & 4 & 0 & 7 & 8 \\
        \end{tabular}
    }
    \]
    This semigroup $S$ is aperiodic, but
    the space $BS$ is a simply connected Moore space
    of type $\mathsf{M}(C_2, 3)$.
    In particular, $H_3(S)=C_2$, and $H_i(S) = 0$ for $i\neq 0, 3$.
    \begin{proof}
        First note that $\pi_1(BS)=GS=G(K(S))=G(\{x_0, x_1, x_2, x_3\})=1$,
        so $BS$ is simply connected. To compute homology,
        the \texttt{fast\_semigroup\_homology} package finds
        the following finite projective resolution of $\Z$ as a $\Z S\one$-module:
        \[
            \Z
            \xleftarrow{\epsilon} \Z S\one x_0
            \xleftarrow{\partial_1} \Z S\one x_5
            \xleftarrow{\partial_2} \begin{matrix}
                \Z S\one x_7 \\
                \oplus \Z S\one
            \end{matrix}
            \xleftarrow{\partial_3} \begin{matrix}
                \Z S\one x_0 \\
                \oplus \Z S\one x_7 \\
                \oplus \Z S\one x_5 \\
                \oplus \Z S\one x_7
            \end{matrix}
            \xleftarrow{\partial_4} \begin{matrix}
                \Z S\one x_0 \\
                \oplus \Z S\one x_0 \\
                \oplus \Z S\one x_0
            \end{matrix}
        \]
        The maps are defined as follows:
        \begin{align*}
            \partial_1(a)&=a(x_0-x_2)\\
            \partial_2(a,b)&=ax_4 + b(x_5-x_6) \\
            \partial_3(a,b,c,d) &= \big(b(x_0-x_7), \, ax_1 + bx_8 + cx_5 + d(x_7 + x_8)\big) \\
            \partial_4(a,b,c) &= \big(2ax_0,\, bx_0,\, bx_0+2cx_0,\, -ax_1-bx_0-cx_0\big).
        \end{align*}
        Removing the augmentation and tensoring with $\Z$ results in the following
        complex, which has the stated homology:
        \[
            \Z \xleftarrow{0} \Z
            \xleftarrow{\left(\begin{smallmatrix}1 & 0\end{smallmatrix}\right)} \Z^2
            \xleftarrow{\left(\begin{smallmatrix}0 & 0 & 0 & 0 \\ 1 & 1 & 1 & 2\end{smallmatrix}\right)} \Z^4
            \xleftarrow{\left(\begin{smallmatrix}2 & 0 & 0 \\ 0 & 1 & 0 \\ 0 & 1 & 2 \\ -1 & -1 & -1\end{smallmatrix}\right)} \Z^3.
            \qedhere
        \]
    \end{proof}
\end{example}

\begin{example}\label{ex:Moore_C2_2}
    Consider the semigroup $S = \{x_0, \dots, x_{11}\}$
    with operation defined by $x_ix_j=x_k$ for every entry $k$ in row $i$ and column $j$
    in the following table\footnote{
        This is a transposed and permuted version of
        the table at index \texttt{297208} of the \texttt{order13.hdf5}
        file in the ``\texttt{...bounded\_qdiag}'' results folder from the \texttt{semisearch} repository,
        with the adjoined identity removed.
        Since the first version of this paper,
        this construction has been generalized by \cite{realizingMooore}.
    }:
    \[
    \text{
    \footnotesize
    \tabcolsep=0.1cm
        \begin{tabular}{c|cccccccccccc}
        $\cdot$  & 0 & 1 & 2 & 3 & 4 & 5 & 6 & 7 & 8 & 9 & 10 & 11 \\ \hline
        0  & 0&0&0&0&4&4&4&4&0&0&4&4 \\
        1  & 1&1&1&1&5&5&5&5&1&1&5&5 \\
        2  & 2&2&2&2&6&6&6&6&2&2&6&6 \\
        3  & 3&3&3&3&7&7&7&7&3&3&7&7 \\
        4  & 0&0&0&0&4&4&4&4&4&4&0&0 \\
        5  & 1&1&1&1&5&5&5&5&5&5&1&1 \\
        6  & 2&2&2&2&6&6&6&6&6&6&2&2 \\
        7  & 3&3&3&3&7&7&7&7&7&7&3&3 \\
        8  & 0&1&2&0&4&5&6&4&8&9&10&11 \\
        9  & 0&1&2&1&4&5&6&5&8&9&10&11 \\
        10 & 1&0&2&1&5&4&6&5&10&11&8&9 \\
        11 & 1&0&2&0&5&4&6&4&10&11&8&9 \\
        \end{tabular}
    }
    \]
    The classifying space $BS$ is a simply connected Moore space
    of type $\mathsf{M}(C_2, 2)$. In particular, $H_2(S)=C_2$ and $H_i(S)=0$
    for $i\neq 0, 2$.
\end{example}
\begin{proof}
    Since $x_0, \dots, x_{11}$ are idempotent, we have
    \[
        \pi_1(BS)=GS=G(K(S)) = G(\{x_0, \dots, x_{11}\}) = 1,
    \]
    so $BS$ is simply connected.
    To compute homology,
    the \texttt{fast\_semigroup\_homology} package finds
    the following finite projective resolution of $\Z$ as a $\Z S\one$-module:
    \[
        \Z
        \xleftarrow{\epsilon} \Z S\one x_0
        \xleftarrow{\partial_1} \Z S\one
        \xleftarrow{\partial_2} \Z S\one x_0 \oplus \Z S\one x_8
        \xleftarrow{\partial_3} \Z S\one x_0.
    \]
    The boundary maps above are:
    \begin{align*}
        \partial_1(a) &= a(x_2-x_3)\\
        \partial_2(a,b) &= ax_4 + b(x_8-x_{11})\\
        \partial_3(a) &= (0, a(x_0 + x_4)).
    \end{align*}
    After removing the augmentation and tensoring with $\Z$, we are left with
    the following complex, which has the stated homology:
    \[
        \Z \xleftarrow{0} \Z
        \xleftarrow{\left(\begin{smallmatrix}1 & 0\end{smallmatrix}\right)} \Z^2
        \xleftarrow{\left(\begin{smallmatrix}0 \\ 2\end{smallmatrix}\right)} \Z.
        \qedhere
    \]
\end{proof}

\begin{example}\label{ex:1494640}
    Consider the semigroup $S = \{x_0, \dots, x_{9}\}$
    with operation defined by $x_ix_j=x_k$ for every entry $k$ in row $i$ and column $j$
    in the following table\footnote{
        This is the table at index \texttt{599445} of the \texttt{order11.hdf5}
        file in the
        ``\texttt{monoids no monoid 1sided ideal by min ideal and diagonal and units}''
        results folder from the \texttt{semisearch} repository,
        with the adjoined identity removed.
    }:
    \[
    \text{
    \footnotesize
    \tabcolsep=0.1cm
        \begin{tabular}{c|cccccccccc}
        $\cdot$  & 0 & 1 & 2 & 3 & 4 & 5 & 6 & 7 & 8 & 9 \\ \hline
        0        & 0 & 1 & 0 & 1 & 0 & 0 & 0 & 0 & 0 & 0 \\
        1        & 0 & 1 & 0 & 1 & 0 & 0 & 0 & 0 & 0 & 0 \\
        2        & 2 & 3 & 2 & 3 & 2 & 2 & 2 & 2 & 2 & 2 \\
        3        & 2 & 3 & 2 & 3 & 2 & 2 & 2 & 2 & 2 & 2 \\
        4        & 0 & 1 & 0 & 1 & 0 & 0 & 0 & 0 & 0 & 0 \\
        5        & 0 & 1 & 0 & 1 & 0 & 0 & 0 & 0 & 0 & 0 \\
        6        & 0 & 1 & 0 & 1 & 0 & 0 & 0 & 0 & 0 & 0 \\
        7        & 0 & 1 & 0 & 1 & 0 & 0 & 0 & 5 & 4 & 6 \\
        8        & 0 & 1 & 0 & 1 & 0 & 0 & 0 & 6 & 4 & 0 \\
        9        & 0 & 1 & 0 & 1 & 0 & 0 & 0 & 6 & 5 & 4 \\
        \end{tabular}
    }
    \]
    Note that this semigroup is a union of a $2$-by-$2$ rectangular band on $\{0,1,2,3\}$
    and a 3-nilpotent semigroup on $\{0,4,5,6,7,8,9\}$.
    It has the following first several homology groups $H_0(S), \dots, H_{10}(S)$:
    \begin{align*}
    &\text{
        \begin{tabular}{c|cccccccc}
            $i$      &    0 & 1 & 2 & 3 & 4 & 5 & 6 & 7\\ \hline
            $H_i(S)$ & $\Z$ & $0$ & $\Z$ & $\Z^3$ & $\Z^6$ & $\Z^9$ & $\Z^{9}\times C_{1494640}$ & $\Z^{27} \times C_{17}$
        \end{tabular}
    }\\
    &\text{
        \begin{tabular}{c|ccc}
            $i$      & 8 & 9 & 10 \\ \hline
            $H_i(S)$ & $\Z^{81}\times C_{17}$ & $\Z^{162} \times C_{17}$ & $\Z^{243} \times C_{1494640} \times C_{17}$\\
        \end{tabular}
    }
    \end{align*}
    We omit the large projective resolution used to calculate the above,
    but the results can be replicated using the command
    ``\texttt{python3 -m fast\_semigroup\_homology}'' with the flags
    ``\texttt{-i "0101000000;[...];0101000654"}'' to input the given multiplication table
    and ``$\texttt{-d 10}$'' to compute up to dimension 10. This
    calculation took approximately 5 minutes to complete.

    The large torsion part
    $C_{1494640} \cong C_{157} \times C_{17} \times C_{7} \times C_{5} \times C_{2^4}$
    notably includes elements of finite order (including prime order)
    much larger than the order of the original semigroup.
\end{example}

\section{Joins and Suspensions}\label{sec:suspensions}
The results of this section are inspired by the following example.
\begin{example}\label{ex:sus_C2}
    Consider the monoid on $M=\{0, \dots, 5\}$
    with the following multiplication table:
    \[
    \text{
    \footnotesize
    \tabcolsep=0.1cm
        \begin{tabular}{c|cccccc}
        $\cdot$  & 0 & 1 & 2 & 3 & 4 & 5 \\ \hline
        0        & 0 & 1 & 0 & 1 & 0 & 1\\
        1        & 0 & 1 & 0 & 1 & 1 & 0\\
        2        & 2 & 3 & 2 & 3 & 2 & 3\\
        3        & 2 & 3 & 2 & 3 & 3 & 2\\
        4        & 0 & 1 & 2 & 3 & 4 & 5\\
        5        & 0 & 1 & 2 & 3 & 5 & 4
        \end{tabular}
    }
    \]
    The \texttt{fast\_semigroup\_homology} package
    computes the following homology for $M$:
    \[
    \text{
        \begin{tabular}{c|ccccccccccc}
            $i$      &    0 & 1 & 2 & 3 & 4 & 5 & 6 & 7 & 8 & 9 & 10\\ \hline
            $H_i(M)$ & $\Z$ & $0$ & $C_2$ & $0$ & $C_2$ & $0$ & $C_2$ & $0$ & $C_2$ & $0$ & $C_2$
        \end{tabular}
    }
    \]
    This section will show that in fact, $BM$ is homotopy equivalent to
    the suspension $\Sigma(\mathbb{R}P^\infty)$ of the infinite real projective space,
    i.e., $BM \simeq \Sigma(BC_2)$. Note that
    because suspensions have trivial cup products
    but $M$ has infinitely many nontrivial (co)homology groups,
    the cohomology ring $H^*(BM)$ cannot be finitely generated as a ring.
\end{example}
\begin{definition}
    For a monoid $M$ and a set $Y$, define a monoid $J^Y(M)$
    on the $(1+|Y|)|M|$ symbols $M \sqcup \{ym : y \in Y, m \in M\}$,
    with operation ``$\cdot$'' defined as follows:
    \begin{align*}
        m \cdot m' &\coloneq mm' \in M & m \cdot y'm' &\coloneq y'm' \\
        ym \cdot m' &\coloneq y(mm')  & ym \cdot y'm' &\coloneq ym'
    \end{align*}
    For $y \in Y$ identify $y \coloneq y1 \in J^Y(M)$, so $Y\subseteq J^Y(M)$
    is a subsemigroup satisfying $yy'=y$.
\end{definition}
Note that $J^Y(M)$ is the union of the monoid $M$
with the $|Y|$-by-$|M|$ rectangular band $YM$,
where $M$ acts as the identity on the left of $YM$
and acts as multiplication in $M$ on the right of $YM$.
This section will prove Theorem~\ref{thm:join_semigroup},
which states that the classifying space $B\left(J^Y(M)\right)$
is homotopy equivalent to the join of $BM$ with the discrete set $Y$.

The \textit{join} of
topological spaces is discussed in Section 5.7 of \cite{brownTopology},
and joins of simplicial sets are discussed in \cite{joins}.
For simplicial sets $X, X'$, each regarded as augmented
with a unique $(-1)$-cell $\varnothing$,
the set of $n$-simplices in the join $X * X'$ is
\[
    (X * X')^n \cong \bigsqcup_{\substack{(a+1)+(b+1)=(n+1)\\ a, b \geq -1}}
    X^a \times X'^b.
\]
Each simplex is denoted as $\sigma^a * \tau^b \in X^a \times X'^b \subseteq (X * X')^{a+b+1}$.
Face maps are given by
\[
    d_i^{a+b+1}(\sigma^a * \tau^b) \coloneq \begin{cases}
        d_i^a\sigma^a * \tau^b & 0 \leq i \leq a \\
        \sigma^a * d^b_{i-(a+1)}\tau^b & a+1\leq i \leq a + b + 1,
    \end{cases}
\]
and degeneracy maps defined similarly. Joins respect geometric realization:
$|X * X'| \cong |X| * |X'|$ \cite[Theorem 4.3]{joins}.

The corollaries from Section~\ref{sec:intro} are immediate
from Theorem~\ref{thm:join_semigroup}:
\begin{proof}[Proof of Corollary~\ref{cor:wedge_of_suspensions}]
    By standard topological properties of joins,
    \[
        BM * Y^{\text{discrete}}
        \simeq \Sigma\!\left(BM \wedge Y^\text{discrete}\right)
        \simeq \Sigma\!\left(\bigvee_{y \in Y \setminus \{y_0\}} BM\right)
        \cong \bigvee_{y \in Y \setminus \{y_0\}} \Sigma(BM).
        \qedhere
    \]
\end{proof}
\begin{proof}[Proof of Corollary~\ref{cor:suspension}]
    Taking $Y=\{1,2\}$, it is immediate from Theorem~\ref{thm:join_semigroup} that
    $B(J^{1,2}(M))\simeq\Sigma(BM)$.
\end{proof}
\begin{proof}[Proof of Corollary~\ref{cor:wedge_of_spheres}]
    To create a finite monoid $M$ with $BM\simeq \bigvee_{k=1}^r \mathbb{S}^n$,
    we can repeatedly apply $J^{\{1,2\}}$ to reduce the problem
    to finding a finite monoid $M$ with $BM\simeq \bigvee_{k=1}^r \mathbb{S}^2$.
    This is supplied by a $(r+1)$-by-$2$ rectangular band with an adjoined unit,
    directly generalizing Example~\ref{ex:rect32}.
    Alternately, we can appeal to \cite{steinberg2024}.
\end{proof}

To prove Theorem~\ref{thm:join_semigroup}, we use the collapsing
scheme defined in \cite{Brown92} for the nerve of a monoid with a
complete rewriting system. See \cite{SOK} for more details
on string rewriting systems for monoids.
\begin{proposition}\label{prop:good_CRS}
    $J^Y(M)$ is a monoid with a complete rewriting system
    on the alphabet $(M\setminus \{1\}) \sqcup Y$ with rewriting rules
\begin{align*}
    mm' &\to m'' && \text{for each $m,m' \in M \setminus \{1\}$ where $mm'=m''$ in $M$},\\
    my &\to y && \text{for each $m \in M\setminus\{1\}, y \in Y$, and}\\
    yy' &\to y && \text{for each $y,y' \in Y$}.
\end{align*}
\end{proposition}
\begin{proof}
    There are no infinite chains of reductions because each rule reduces
    the length of a string. The system is confluent because
    the critical pairs converge:
    \begin{itemize}
        \item $(mm')y\to m''y \to y$ agrees with $m(m'y)\to my \to y$.
        \item $(my)y' \to yy' \to y$ agrees with $m(yy')\to my\to y$.
        \item Critical pairs only in $S$ or only in $Y$ certainly converge.
    \end{itemize}
    The irreducible words in this system take the form of $m$, $y$, $ym$ or $1$
    for any $m \in M\setminus \{1\}$ and $y \in Y$,
    and the rewriting rules give exactly the operation on $J^Y(M)$.
\end{proof}
\begin{lemma}\label{lem:brown_Q}
    There is a homotopy equivalence $B(J^Y(M)) \to Q$,
    where $Q$ is the subcomplex of $B(J^Y(M))$ consisting
    of cells of the form $[m_1|\cdots|m_k|y_1| \cdots| y_\ell]$,
    where each $m_i \in M\setminus\{1\}$ and $y_i \in Y$, and $k, \ell \geq 0$.
\end{lemma}
\begin{proof}
    Using the complete rewriting system of Proposition~\ref{prop:good_CRS},
    the collapsing scheme of \cite{Brown92}
    defines a homotopy equivalence from $B(J^Y(M))$
    onto a certain quotient complex with a cell for each \textit{essential}
    cell of $B(J^Y(M))$.
    Using our notation from Section~\ref{sub:classifying},
    Brown's essential cells
    are lists $[w_1|\cdots|w_n]$ of irreducible words over the alphabet,
    with the first entry $w_1$ a single letter, the concatenation
    $w_iw_{i+1}$ of any pair of adjacent entries reducible,
    and no proper prefix of any $w_iw_{i+1}$ reducible.

    For our system, the essential cells cannot include $[\cdots|m|ym'| \cdots]$
    and cannot include $[\cdots|y|y'm'|\cdots]$ because $my$ and $yy'$
    are reducible proper prefixes of the concatenations.
    Thus, because these longer irreducible words cannot follow a single-letter word,
    all entries of our essential cells are single letters.
    Because any $ym$ is irreducible,
    our essential cells cannot include $[\cdots|y|m|\cdots]$ either,
    so we are left with only the listed cells of $Q$.

    For an arbitrary complete rewriting system,
    Brown's quotient complex is a more general CW complex.
    However, because our essential cells and their degeneracies
    form a simplicial set subcomplex of $B(J^Y(M))$,
    no stage of the collapsing scheme
    alters the face maps on $Q$, so we can identify $Q$
    with the stated simplicial set subcomplex.
\end{proof}

\begin{proof}[Proof of Theorem~\ref{thm:join_semigroup}]
    Using the simplicial set $Q$ from Lemma~\ref{lem:brown_Q},
    we show the following zigzag of homotopy equivalences:
    \[
        BM * Y^\text{discrete}
        \underset{(1)}{\leftarrow}
        BM * \bigsqcup_{y \in Y} {BY}
        \underset{(2)}{\rightarrow}
        \frac{BM * \bigsqcup_{y \in Y} {BY}}{() * \bigsqcup_{y \in Y} {BY}}
        \underset{(3)}{\cong}
        \frac{Q}{BY}
        \underset{(4)}{\leftarrow}
        Q
        \underset{(5)}{\leftarrow}
        B\!\left(J^Y(M)\right).
    \]
    Here, $Y \subseteq J^Y(M)$ has operation $yy' = y$,
    so $BY$ is contractible by Corollary~\ref{cor:zero_contract},
    and so (4) is an equivalence.
    Joins respect equivalences,
    so (1) is an equivalence.
    A join with a single point $()$ is
    a cone, and hence is contractible, so (2) is an equivalence.
    The map (5) is an equivalence by the lemma,
    so it remains to give the isomorphism (3) of simplicial sets.
    The correspondence (3) on nondegenerate cells is
    \begin{align*}
        \text{the unique $0$-cell}
        &\longleftrightarrow \text{the unique $0$-cell}
        \\
        \big[[m_1|\cdots|m_k] * \varnothing\big]
        &\longleftrightarrow
        \big[[m_1|\cdots| m_k]\big]
        &&\text{ for $k \geq 1$}
        \\
        \big[[m_1|\cdots|m_k] * \iota_{y_0} [y_1|\cdots|y_\l]\big]
        &\longleftrightarrow \big[[m_1|\cdots|m_k|y_0|y_1|\cdots|y_\l]\big]
        &&\text{ for $k \geq 1, \l\geq 0$}.
    \end{align*}
    For the last case, the face map $d_k$ is respected because
    $m_ky_0 = y_0$, and
    the face map $d_{k+1}$
    is respected because $y_0y_1 = y_0$.
\end{proof}

\printbibliography

\begin{table}[p]
    \caption{The format for multiplication tables of monoids
    with minimal ideal a 2-by-2 rectangular band
    on $\{0,1,2,3\}$, and with group of units a cyclic group of order 2 on the elements
    $\{6, 7\}$. Vertical bars separate possibilities for each entry.
    The possibilities listed define the minimal left ideals $\{0,2\}$, $\{1,3\}$
    and minimal right ideals $\{0,1\}$, $\{2,3\}$
    and ensure that the non-units form an ideal.
    Additional constraints can enforce associativity, fix specific diagonal entries,
    and ensure each monoid isomorphism type is represented only by its
    lexicographically minimal table of this form.}
    \label{tab:KS_is_rect22}
    \begin{tabular}{c|cccccccc}
    $\cdot$ & 0 & 1 & 2 & 3 & 4 & 5 & 6 & 7\\ \hline
    0     &0&1&0&1&$0|1$&$0|1$&0&$0|1$\\
    1     &0&1&0&1&$0|1$&$0|1$&1&$0|1$\\
    2     &2&3&2&3&$2|3$&$2|3$&2&$2|3$\\
    3     &2&3&2&3&$2|3$&$2|3$&3&$2|3$\\
    4     &$0|2$&$1|3$&$0|2$&$1|3$&$0|1|2|3|4|5$&$0|1|2|3|4|5$&4&$0|1|2|3|4|5$\\
    5     &$0|2$&$1|3$&$0|2$&$1|3$&$0|1|2|3|4|5$&$0|1|2|3|4|5$&5&$0|1|2|3|4|5$\\
    6     &0&1&2&3&4&5&6&7\\
    7     &$0|2$&$1|3$&$0|2$&$1|3$&$0|1|2|3|4|5$&$0|1|2|3|4|5$&7&6\\
    \end{tabular}
\end{table}

\begin{table}[p]
    \caption{
    For each order $\leq 12$, listed are the number of semigroups of order $n$,
    the number of monoids of order $n$, and the number of monoids of order $n$
    with no idempotent $e\neq 1$ such that $eMe=eM$ or $eMe=Me$,
    so as not to satisfy the hypothesis of Theorem~\ref{thm:eSeEquivalent}.
    All of these are counted up to isomorphism and anti-isomorphism.
    The semigroups data is from \cite{smallsemi,order10},
    with the lower bounds provided by 3-nilpotent semigroups \cite[]{nilpotents}.
    The monoids data is from \cite{mon8910}, with lower bounds
    provided by the number of semigroups one order smaller.
    The last column's data is from \texttt{semisearch}.
    }
    \label{tab:counts_table}
    \centering
    \begin{tabular}{rrrr}
        \toprule
        Order & Semigroups & Monoids & Monoids, no $eMe$ \\
        \midrule
        1 & 1        & 1        & 1 \\
        2 & 4        & 2        & 1 \\
        3 & 18       & 6        & 1 \\
        4 & 126      & 27       & 2 \\
        5 & 1\,160   & 156      & 2 \\
        6 & 15\,973  & 1\,373   & 7 \\
        7 & 836\,021 & 17\,730  & 25 \\
        8 & 1\,843\,120\,128        & 858\,977         & 256 \\
        9 & 52\,989\,400\,714\,478  & 1\,844\,075\,697 & 3\,665 \\
        10& 12\,418\,001\,077\,381\,302\,684
            & 52\,991\,253\,973\,742
            & 71\,916 \\
        11& $>2.6\times 10^{25}$ & $> 1.2\times 10^{19}$ & 2\,232\,321 \\
        12& $>1.0\times 10^{33}$ & $> 2.6\times 10^{25}$ & 1\,974\,639\,821 \\
        \bottomrule
    \end{tabular}
\end{table}

\begin{table}[p]
    \caption{
    Listed are six integral homology groups
    for each semigroup of order at most $6$.
    A dot denotes the trivial group, and $C_k$ denotes the cyclic group of order $k$.
    Within the table for each order,
    each distinct list of homology groups is represented only once,
    and the number of semigroups up to isomorphism and anti-isomorphism
    with the given homology groups is recorded in the column labeled ``$\#$''.
    Further results are available in the ``results'' folder of \cite{fast_semigroup_homology}.
    The large fraction of semigroups with the same homology
    as some group provides motivation for restricting calculations to
    only those monoids for which Theorem~\ref{thm:eSeEquivalent} does not apply,
    as described in Section~\ref{sec:finding}.
    }
    \label{tab:semigroup_homology_tables}
\setlength{\columnsep}{0cm}
\begin{multicols}{2}
    \begin{tabular}{|cccccc|c|}
        \multicolumn{7}{l}{\textbf{Semigroups of Order 0:}\rule{0pt}{2.6ex}} \\
        \hline
        $H_1$ & $H_2$ & $H_3$ & $H_4$ & $H_5$ & $H_6$ & $\#$ \\
        \hline
        $\cdot$ & $\cdot$ & $\cdot$ & $\cdot$ & $\cdot$ & $\cdot$ & 1 \\
        \hline
    \end{tabular}
    \begin{tabular}{|cccccc|c|}
        \multicolumn{7}{l}{\textbf{Semigroups of Order 1:}\rule{0pt}{2.6ex}} \\
        \hline
        $H_1$ & $H_2$ & $H_3$ & $H_4$ & $H_5$ & $H_6$ & $\#$ \\
        \hline
        $\cdot$ & $\cdot$ & $\cdot$ & $\cdot$ & $\cdot$ & $\cdot$ & 1 \\
        \hline
    \end{tabular}
    \begin{tabular}{|cccccc|c|}
        \multicolumn{7}{l}{\textbf{Semigroups of Order 2:}\rule{0pt}{2.6ex}} \\
        \hline
        $H_1$ & $H_2$ & $H_3$ & $H_4$ & $H_5$ & $H_6$ & $\#$ \\
        \hline
        $\cdot$ & $\cdot$ & $\cdot$ & $\cdot$ & $\cdot$ & $\cdot$ & 3 \\
        $C_2$ & $\cdot$ & $C_2$ & $\cdot$ & $C_2$ & $\cdot$ & 1 \\
        \hline
    \end{tabular}
    \begin{tabular}{|cccccc|c|}
        \multicolumn{7}{l}{\textbf{Semigroups of Order 3:}\rule{0pt}{2.6ex}} \\
        \hline
        $H_1$ & $H_2$ & $H_3$ & $H_4$ & $H_5$ & $H_6$ & $\#$ \\
        \hline
        $\cdot$ & $\cdot$ & $\cdot$ & $\cdot$ & $\cdot$ & $\cdot$ & 14 \\
        $C_2$ & $\cdot$ & $C_2$ & $\cdot$ & $C_2$ & $\cdot$ & 3 \\
        $C_3$ & $\cdot$ & $C_3$ & $\cdot$ & $C_3$ & $\cdot$ & 1 \\
        \hline
    \end{tabular}
    \begin{tabular}{|cccccc|c|}
        \multicolumn{7}{l}{\textbf{Semigroups of Order 4:}\rule{0pt}{2.6ex}} \\
        \hline
        $H_1$ & $H_2$ & $H_3$ & $H_4$ & $H_5$ & $H_6$ & $\#$ \\
        \hline
        $\cdot$ & $\cdot$ & $\cdot$ & $\cdot$ & $\cdot$ & $\cdot$ & 102 \\
        $C_2$ & $\cdot$ & $C_2$ & $\cdot$ & $C_2$ & $\cdot$ & 18 \\
        $C_3$ & $\cdot$ & $C_3$ & $\cdot$ & $C_3$ & $\cdot$ & 3 \\
        $\cdot$ & $\Z$ & $\cdot$ & $\cdot$ & $\cdot$ & $\cdot$ & 1 \\
        $C_2^{2}$ & $C_2$ & $C_2^{3}$ & $C_2^{2}$ & $C_2^{4}$ & $C_2^{3}$ & 1 \\
        $C_4$ & $\cdot$ & $C_4$ & $\cdot$ & $C_4$ & $\cdot$ & 1 \\
        \hline
    \end{tabular}
    \begin{tabular}{|cccccc|c|}
        \multicolumn{7}{l}{\textbf{Semigroups of Order 5:}\rule{0pt}{2.6ex}} \\
        \hline
        $H_1$ & $H_2$ & $H_3$ & $H_4$ & $H_5$ & $H_6$ & $\#$ \\
        \hline
        $\cdot$ & $\cdot$ & $\cdot$ & $\cdot$ & $\cdot$ & $\cdot$ & 996 \\
        $C_2$ & $\cdot$ & $C_2$ & $\cdot$ & $C_2$ & $\cdot$ & 136 \\
        $C_3$ & $\cdot$ & $C_3$ & $\cdot$ & $C_3$ & $\cdot$ & 17 \\
        $C_4$ & $\cdot$ & $C_4$ & $\cdot$ & $C_4$ & $\cdot$ & 4 \\
        $C_2^{2}$ & $C_2$ & $C_2^{3}$ & $C_2^{2}$ & $C_2^{4}$ & $C_2^{3}$ & 3 \\
        $\cdot$ & $\Z$ & $\cdot$ & $\cdot$ & $\cdot$ & $\cdot$ & 2 \\
        $\cdot$ & $\Z$ & $\Z$ & $\Z$ & $\Z$ & $\Z$ & 1 \\
        $C_5$ & $\cdot$ & $C_5$ & $\cdot$ & $C_5$ & $\cdot$ & 1 \\
        \hline
    \end{tabular}
    \begin{tabular}{|cccccc|c|}
        \multicolumn{7}{l}{\textbf{Semigroups of Order 6:}\rule{0pt}{2.6ex}} \\
        \hline
        $H_1$ & $H_2$ & $H_3$ & $H_4$ & $H_5$ & $H_6$ & $\#$ \\
        \hline
        $\cdot$ & $\cdot$ & $\cdot$ & $\cdot$ & $\cdot$ & $\cdot$ & 14334 \\
        $C_2$ & $\cdot$ & $C_2$ & $\cdot$ & $C_2$ & $\cdot$ & 1420 \\
        $C_3$ & $\cdot$ & $C_3$ & $\cdot$ & $C_3$ & $\cdot$ & 139 \\
        $C_4$ & $\cdot$ & $C_4$ & $\cdot$ & $C_4$ & $\cdot$ & 25 \\
        $C_2^{2}$ & $C_2$ & $C_2^{3}$ & $C_2^{2}$ & $C_2^{4}$ & $C_2^{3}$ & 18 \\
        $\cdot$ & $\Z$ & $\cdot$ & $\cdot$ & $\cdot$ & $\cdot$ & 17 \\
        $\cdot$ & $\Z$ & $\Z$ & $\Z$ & $\Z$ & $\Z$ & 6 \\
        $C_5$ & $\cdot$ & $C_5$ & $\cdot$ & $C_5$ & $\cdot$ & 3 \\
        $\cdot$ & $\Z$ & $\Z^{2}$ & $\Z^{4}$ & $\Z^8$ & $\Z^{16}$ & 3 \\
        $\cdot$ & $\cdot$ & $\Z$ & $\cdot$ & $\cdot$ & $\cdot$ & 2 \\
        $\cdot$ & $\Z$ & $C_2$ & $\cdot$ & $C_2$ & $\cdot$ & 2 \\
        $\cdot$ & $C_2$ & $\cdot$ & $C_2$ & $\cdot$ & $C_2$ & 1 \\
        $\cdot$ & $\Z^{2}$ & $\cdot$ & $\cdot$ & $\cdot$ & $\cdot$ & 1 \\
        $C_2$ & $\cdot$ & $C_6$ & $\cdot$ & $C_2$ & $\cdot$ & 1 \\
        $C_6$ & $\cdot$ & $C_6$ & $\cdot$ & $C_6$ & $\cdot$ & 1 \\
        \hline
    \end{tabular}
\end{multicols}
\end{table}

\begin{table}[p]
    \caption{
    Listed are six integral homology groups
    of those monoids of orders at most $7$
    with no idempotent $e\neq 1$ such that $eMe=eM$ or $eMe=Me$,
    so as not to satisfy the hypothesis of Theorem~\ref{thm:eSeEquivalent}.
    The format follows Table~\ref{tab:semigroup_homology_tables}.
    Note that restricting from the arbitrary semigroups
    of Table~\ref{tab:semigroup_homology_tables}
    to only this restricted class of monoids
    greatly reduces duplication of homotopy types.
    Again, further results are available in the ``results'' folder
    of \cite{fast_semigroup_homology}.
    }
    \label{tab:monoid_homology_tables}
\setlength{\columnsep}{0cm}
\begin{multicols}{2}
    \begin{tabular}{|cccccc|c|}
        \multicolumn{7}{l}{\textbf{Monoids of Order 1, no $eMe$:}\rule{0pt}{2.6ex}} \\
        \hline
        $H_1$ & $H_2$ & $H_3$ & $H_4$ & $H_5$ & $H_6$ & $\#$ \\
        \hline
        $\cdot$ & $\cdot$ & $\cdot$ & $\cdot$ & $\cdot$ & $\cdot$ & 1 \\
        \hline
    \end{tabular}
    \begin{tabular}{|cccccc|c|}
        \multicolumn{7}{l}{\textbf{Monoids of Order 2, no $eMe$:}\rule{0pt}{2.6ex}} \\
        \hline
        $H_1$ & $H_2$ & $H_3$ & $H_4$ & $H_5$ & $H_6$ & $\#$ \\
        \hline
        $C_2$ & $\cdot$ & $C_2$ & $\cdot$ & $C_2$ & $\cdot$ & 1 \\
        \hline
    \end{tabular}
    \begin{tabular}{|cccccc|c|}
        \multicolumn{7}{l}{\textbf{Monoids of Order 3, no $eMe$:}\rule{0pt}{2.6ex}} \\
        \hline
        $H_1$ & $H_2$ & $H_3$ & $H_4$ & $H_5$ & $H_6$ & $\#$ \\
        \hline
        $C_3$ & $\cdot$ & $C_3$ & $\cdot$ & $C_3$ & $\cdot$ & 1 \\
        \hline
    \end{tabular}
    \begin{tabular}{|cccccc|c|}
        \multicolumn{7}{l}{\textbf{Monoids of Order 4, no $eMe$:}\rule{0pt}{2.6ex}} \\
        \hline
        $H_1$ & $H_2$ & $H_3$ & $H_4$ & $H_5$ & $H_6$ & $\#$ \\
        \hline
        $C_2^{2}$ & $C_2$ & $C_2^{3}$ & $C_2^{2}$ & $C_2^{4}$ & $C_2^{3}$ & 1 \\
        $C_4$ & $\cdot$ & $C_4$ & $\cdot$ & $C_4$ & $\cdot$ & 1 \\
        \hline
    \end{tabular}
    \begin{tabular}{|cccccc|c|}
        \multicolumn{7}{l}{\textbf{Monoids of Order 5, no $eMe$:}\rule{0pt}{2.6ex}} \\
        \hline
        $H_1$ & $H_2$ & $H_3$ & $H_4$ & $H_5$ & $H_6$ & $\#$ \\
        \hline
        $\cdot$ & $\Z$ & $\cdot$ & $\cdot$ & $\cdot$ & $\cdot$ & 1 \\
        $C_5$ & $\cdot$ & $C_5$ & $\cdot$ & $C_5$ & $\cdot$ & 1 \\
        \hline
    \end{tabular}
    \begin{tabular}{|cccccc|c|}
        \multicolumn{7}{l}{\textbf{Monoids of Order 6, no $eMe$:}\rule{0pt}{2.6ex}} \\
        \hline
        $H_1$ & $H_2$ & $H_3$ & $H_4$ & $H_5$ & $H_6$ & $\#$ \\
        \hline
        $\cdot$ & $\Z$ & $\cdot$ & $\cdot$ & $\cdot$ & $\cdot$ & 1 \\
        $\cdot$ & $\Z$ & $\Z$ & $\Z$ & $\Z$ & $\Z$ & 1 \\
        $\cdot$ & $\Z$ & $C_2$ & $\cdot$ & $C_2$ & $\cdot$ & 2 \\
        $\cdot$ & $C_2$ & $\cdot$ & $C_2$ & $\cdot$ & $C_2$ & 1 \\
        $C_2$ & $\cdot$ & $C_6$ & $\cdot$ & $C_2$ & $\cdot$ & 1 \\
        $C_6$ & $\cdot$ & $C_6$ & $\cdot$ & $C_6$ & $\cdot$ & 1 \\
        \hline
    \end{tabular}
{
\small
    \begin{tabular}{|cccccc|c|}
        \multicolumn{7}{l}{\textbf{Monoids of Order 7, no $eMe$:}\rule{0pt}{2.6ex}} \\
        \hline
        $H_1$ & $H_2$ & $H_3$ & $H_4$ & $H_5$ & $H_6$ & $\#$ \\
        \hline
        $\cdot$ & $\cdot$ & $\Z$ & $\cdot$ & $\cdot$ & $\cdot$ & 2 \\
        $\cdot$ & $\Z$ & $\cdot$ & $\cdot$ & $\cdot$ & $\cdot$ & 10 \\
        $\cdot$ & $\Z$ & $\Z$ & $\Z$ & $\Z$ & $\Z$ & 5 \\
        $\cdot$ & $\Z$ & \scalebox{0.6}{$\Z \times C_2$} & \scalebox{0.6}{$\Z \times C_2^2$} & \scalebox{0.6}{$\Z \times C_2^5$} & \scalebox{0.6}{$\Z \times C_2^{10}$} & 1\\
        $\cdot$ & $\Z$ & $\Z^{2}$ & $\Z^{4}$ & $\Z^8$ & $\Z^{16}$ & 3 \\
        $\cdot$ & $\Z$ & $C_2$ & $\cdot$ & $C_2$ & $\cdot$ & 1 \\
        $\cdot$ & $\Z$ & $C_3$ & $\cdot$ & $C_3$ & $\cdot$ & 1 \\
        $\cdot$ & $\Z^{2}$ & $\cdot$ & $\cdot$ & $\cdot$ & $\cdot$ & 1 \\
        $C_7$ & $\cdot$ & $C_7$ & $\cdot$ & $C_7$ & $\cdot$ & 1 \\
        \hline
    \end{tabular}
}
\end{multicols}
\end{table}

\begin{amssidewaystable}
    \centering
    \caption{
        Benchmarks for various methods of computing homology of finite monoids.
        Each row $n$ shows
        the time in milliseconds to compute the homology groups $H_1, \dots, H_n$
        of 39 specific monoids of orders $\leq 7$.
        Columns under \texttt{fast\_semigroup\_homology} differ only
        by how kernels were computed, with ``Default'' using \texttt{mutable\_lattice}.
        The HAP columns use \cite{HAP},
        and the SageMath column directly constructs the Delta set $B_\Sgrp M$
        from Section~\ref{sub:classifying}.
        A ``$\boxtimes$''
        indicates that the calculation did not finish
        after 10 minutes.
        See the \texttt{benchmarks} folder of \texttt{fast\_semigroup\_homology}
        for the benchmarking code. Results were measured on a laptop with
        an Intel i7-8750H CPU, running SageMath 10.5.beta3 on Ubuntu 22.04.3
        with Windows Subsystem for Linux.
    }\label{tab:benchmarks}
    \begin{tabular}{rrrrrrrrrrrrr}
        \toprule
        {}
        & \multicolumn{9}{c}{\texttt{fast\_semigroup\_homology}}
        & \multicolumn{2}{c}{{}}
        & \multicolumn{1}{c}{{}}
        \\
        \cmidrule{2-10}
        {} & {}
        & \multicolumn{4}{c}{Kernels using PARI}
        & \multicolumn{4}{c}{Kernels using SageMath}
        & \multicolumn{2}{c}{HAP}
        & \multicolumn{1}{c}{SageMath}
        \\
        \cmidrule(lr){3-6}
        \cmidrule(lr){7-10}
        \cmidrule(lr){11-12}
        \cmidrule(lr){13-13}
        $n$ & Default
        & {HNF5} & {HNF1}
        & {HNF4} & {MatKerInt}
        & {pari} & {default}
        & {flint} &{padic}
        & {Contracted} & {BarCx}
        & {DeltaCx}
        \\ \midrule
        1 & 10 & 14 & 13 & 14 & 14 & 24 & 30 & 27 & 63 & 1951 & 1819 & 411 \\
        2 & 14 & 22 & 22 & 23 & 23 & 38 & 46 & 54 & 108 & 2287 & 2260 & 11755 \\
        3 & 17 & 30 & 29 & 31 & 32 & 49 & 59 & 74 & 143 & 4346 & 9434 & $\boxtimes$ \\
        4 & 19 & 34 & 34 & 37 & 36 & 56 & 68 & 87 & 167 & 16615 & 553640 &  \\
        5 & 21 & 38 & 38 & 42 & 43 & 63 & 76 & 105 & 193 & 105901 & $\boxtimes$ &  \\
        6 & 23 & 45 & 45 & 54 & 54 & 76 & 89 & 132 & 231 & $\boxtimes$ &  &  \\
        7 & 25 & 61 & 62 & 85 & 85 & 110 & 161 & 203 & 325 &  &  &  \\
        8 & 30 & 94 & 96 & 179 & 187 & 221 & 365 & 419 & 554 &  &  &  \\
        9 & 39 & 158 & 195 & 393 & 412 & 466 & 755 & 824 & 983 &  &  &  \\
        10 & 60 & 299 & 359 & 1240 & 1312 & 1392 & 4654 & 4421 & 4448 &  &  &  \\
        11 & 131 & 680 & 840 & 5378 & 5479 & 5784 & 30159 & 30672 & 30471 &  &  &  \\
        12 & 309 & 1563 & 1858 & 16915 & 17237 & 18202 & 106359 & 106621 & 108978 &  &  &  \\
        13 & 940 & 3867 & 5142 & 67781 & 69795 & 71445 & 554380 & 545685 & 560411 &  &  &  \\
        14 & 4113 & 11418 & 14916 & 356443 & 367500 & 366476 & $\boxtimes$ & $\boxtimes$ & $\boxtimes$ &  &  &  \\
        15 & 14255 & 32500 & 43561 & $\boxtimes$ & $\boxtimes$ & $\boxtimes$ &  &  &  &  &  &  \\
        16 & 33156 & 75946 & 126006 &  &  &  &  &  &  &  &  &  \\
        17 & 62327 & 254586 & 381806 &  &  &  &  &  &  &  &  &  \\
        18 & 106241 & 526970 & $\boxtimes$ &  &  &  &  &  &  &  &  &  \\
        19 & 405994 & $\boxtimes$ &  &  &  &  &  &  &  &  &  &  \\
        20 & $\boxtimes$ &  &  &  &  &  &  &  &  &  &  &  \\
        \bottomrule
    \end{tabular}
\end{amssidewaystable}

\end{document}